\documentstyle[amsfonts,12pt,thmsa,sw20jart,epsf]{article}

\input tcilatex
\QQQ{Language}{
American English
}

\begin{document}

\author{Marco Mackaay \\
Sector de matematica UCEH\\
Universidade do Algarve\\
8000 Faro \\
Portugal\\
e-mail: mmackaay@ualg.pt}
\title{Enhancing an R-matrix. }
\date{27-11-1996}
\maketitle

\begin{abstract}
\noindent In order to construct a representation of the tangle category one
needs an enhanced R-matrix. In this paper we define a sufficient and
necessary condition for enhancement that can be checked easily for any
R-matrix. If the R-matrix can be enhanced, we also show how to construct the
additional data that define the enhancement. As a direct consequence we find
a sufficient condition for the construction of a knot invariant.

\smallskip\ 

\noindent {\sl AMS Subject Classification: }16W30, 57M25.

\noindent {\sl Keywords }${\frak \&}$ {\sl Phrases:} Hopf algebra, quantum
group, enhanced R-matrix, knot invariant, category of tangles.
\end{abstract}

\begin{section}{Introduction}

\noindent The title of this article contains two words that must be
explained. The first one is the word {\it R-matrix}. An R-matrix is a $%
n^2\times n^2$ matrix $R$ that satisfies the Quantum Yang-Baxter equation 
\begin{equation}
R_{12}R_{13}R_{23}=R_{23}R_{13}R_{12}.  \label{YB}
\end{equation}

\noindent Here we have $R_{12}=R\otimes {\rm id}$ and $R_{23}={\rm id}%
\otimes R$. By $R_{13}$ we mean the following. Let V be a vector space of
dimension $n$, with basis $e_1,\ldots ,e_n.$ Then $R:V\otimes V\rightarrow
V\otimes V$ is given by 
\[
R\left( e_i\otimes e_j\right) =R_{ij}^{kl}e_k\otimes e_l. 
\]

\noindent and $R_{13}:V\otimes V\otimes V\rightarrow V\otimes V\otimes V$ by 
\[
R_{13}\left( e_i\otimes e_j\otimes e_k\right) =R_{ik}^{mn}e_m\otimes
e_j\otimes e_n. 
\]

\noindent The complexity of equation $\left( \ref{YB}\right) $ one better
understands when one writes down the QYB-equation in terms of the matrix
entries $\left( R_{cd}^{ab}\right) .$ Equation $\left( \ref{YB}\right) $
becomes 
\begin{equation}
R_{k_1k_2}^{ab}R_{uk_3}^{k_1c}R_{vw}^{k_2k_3}=R_{l_1l_2}^{bc}R_{l_3w}^{al_2}R_{uv}^{l_3l_1}.
\label{YB2}
\end{equation}

\noindent The number of equations in $\left( \ref{YB2}\right) $ is $n^6$,
which makes a classification of all solutions of the QYB-equation very
difficult, if not impossible. Nonetheless there are by now a great number of
R-matrices known. The most famous ones come from the evaluation of the
universal R-matrix of the universal enveloping algebras of the classical
semi-simple complex Lie algebras in their fundamental representations.
Inspired by the R-matrix associated to the Lie algebras $sl(n)$, Hazewinkel 
\cite{HAZ} classified all solutions of $\left( \ref{YB2}\right) $ under the
restriction 
\begin{equation}
R_{cd}^{ab}\neq 0\Longrightarrow \left\{ a,b\right\} =\left\{ c,d\right\} .
\label{HA}
\end{equation}

\noindent The question now arises what to do with these R-matrices. One
thing one can do with them is try to construct knot and link invariants.
Unfortunately we need some additional data for such a construction, which
brings us to the other word in the title that must be explained.

\noindent Turaev \cite{TUR1} showed what extra data are needed for the
construction of a knot invariant. He used the term {\it enhanced R-matrix}.

\begin{definition}
\cite{TUR1}\label{e1} An enhanced R-matrix is a quadruple $\left( S,\mu
,\alpha ,\beta \right) $ consisting of an invertible $n^2\times n^2$ matrix $S$, a $%
n\times n$ matrix $\mu $ and two complex numbers $\alpha ,\beta \in {\bf 
{\Bbb C}}^{*}$ satisfying the following conditions: 
\begin{equation}
S_{12}S_{23}S_{12}=S_{23}S_{12}S_{23},  \label{YB3}
\end{equation}
\begin{equation}
S\left( \mu \otimes \mu \right) =\left( \mu \otimes \mu \right) S,
\label{ENH1}
\end{equation}
\begin{equation}
{\rm Tr}_2\left( S^{\pm 1}\left( \mu \otimes \mu \right) \right) =\alpha
^{\pm 1}\beta \mu .  \label{ENH2}
\end{equation}
\end{definition}

\noindent The solutions of $\left( \ref{YB3}\right) $ have a simple relation
with the solutions of $\left( \ref{YB}\right) $. If $R$ is a solution of $%
\left( \ref{YB}\right) $, then both $PR$ and $RP$ are solutions of $\left( 
\ref{YB3}\right) $, where $P$ is the permutation matrix $P_{cd}^{ab}=\delta
_d^a\delta _c^b$. When written down in matrix entries we get 
\begin{equation}
\left( PR\right) _{cd}^{ab}=R_{cd}^{ba},\bigskip\ \left( RP\right)
_{cd}^{ab}=R_{dc}^{ab}.  \label{S}
\end{equation}

\noindent We say that an R-matrix $R$ can be enhanced if there exists an enhanced 
R-matrix with $S=PR$ or $S=RP$. The symbol ${\rm Tr}_2$ stands for the second trace, which in
terms of matrix entries is defined by 
\begin{equation}
{\rm Tr}_2(A)_c^a=A_{cd}^{ad}.  \label{TR2}
\end{equation}

\noindent This definition is independent of the basis with respect to which $%
A$ is written (see \cite{KAS}). When $\mu $ is invertible, then $\left( \ref
{ENH2}\right) $ is equivalent to 
\begin{equation}
{\rm Tr}_2\left( S^{\pm 1}\left( I\otimes \mu \right) \right) =\alpha
^{\pm }\beta I.  \label{ENH3}
\end{equation}

\noindent If $\left( S,\mu ,\alpha ,\beta \right) $ is an enhanced R-matrix,
then $\left( \alpha ^{-1}S,\beta ^{-1}\mu ,1,1\right) $ is one also. So both
the factors $\alpha $ and $\beta $ can be normalized to 1. Given such an
enhanced R-matrix, Turaev constructed the link invariant 
\begin{equation}  \label{TU}
T_S\left( \xi \right) =\alpha ^{-w\left( \xi \right) }\beta ^{-m}{\rm Tr}%
\left( \rho _S\left( \xi \right) \circ \mu ^{\otimes m}\right) .
\end{equation}

\noindent 
Here $\xi $ is a braid with m strands, $w\left( \xi \right) =\sum
\varepsilon _i$ if $\xi =\sigma _{i_1}^{\varepsilon _1}\cdots \sigma
_{i_r}^{\varepsilon _r}$, where the $\sigma _i$ are the standard generators
of the braid group on m letters $B_m$ and $\rho _S$ is the braid
representation in $\left( {\bf {\Bbb C}}^n\right) ^{\otimes m}$ defined by 
\[
\rho _S\left( \xi \right) =S_{i_1i_1+1}^{\varepsilon _1}\cdots
S_{i_ri_r+1}^{\varepsilon _r}. 
\]

\noindent The invariant is well defined on links because the trace is
actually a Markov trace, which means that its value is independent of the
braid presentation of the link.

In \cite{TUR3} Turaev gives a slightly more restrictive definition
of an enhanced R-matrix in order to get a sufficient condition for the
construction of a representation of the category of tangles. There his definition
is the following:

\begin{definition}
\label{e2}Let $V$ be a finite-dimensional vector space. An enhanced R-matrix
is a pair $\left( S,\mu \right) $ where $S$ is an automorphism of $V\otimes V
$ and $\mu $ an automorphism of $V$ satisfying conditions $\left( \ref{YB3}%
\right) ,$ $\left( \ref{ENH1}\right) $ and $\left( \ref{ENH3}\right) $ with $\alpha=\beta=1$,  
and
additionally 
\begin{equation}
\left( PS^{\mp 1}\right) ^{t_1}\left( I_{V^{*}}\otimes \mu \right) \left(
S^{\pm 1}P\right) ^{t_1}\left( I_{V^{*}}\otimes \mu ^{-1}\right)
=I_{V^{*}\otimes V}.  \label{ENH4}
\end{equation}
\end{definition}

\noindent Here $t_1$ stands for the first transpose, which on matrices is
defined by 
\[
\left( A^{t_1}\right) _{cd}^{ab}=A_{ad}^{cb}. 
\]
With this definition he constructs a unique functor from the category of
tangles to the category of vector spaces such that $F\left( +\right) =V$ and 
$F\left( -\right) =V^{*},$ for every enhanced R-matrix $\left( S,\mu \right) 
$. When restricted to links Turaev's functor gives exactly $\left(\ref{TU}\right).$
We postpone the definition of Turaev's functor to section 4, because it requires some 
of the basic results about the tangle category, which we explain in section 2. 

The question that one asks naturally, after learning the meaning of the
words enhanced R-matrix, is where this enhancement comes from. Turaev \cite
{TUR1} enhanced the R-matrices associated to the classical semi-simple
complex Lie algebras so that they match his definition, but he gave no
explanation for his enhancement. Hazewinkel \cite{HAZ} gave a criterion for
the enhancement of the R-matrices that satisfy the restriction $\left( \ref
{HA}\right) $ so that they match Turaev's definition, but he got his
criterion in a purely combinatorial way that does not reveal the nature of
enhancement. In this paper we explain enhancement in the language of
category theory, which has become the most succesful way of looking at knot
invariants coming from Lie algebras and their R-matrices, and give a simple
criterion for whether an R-matrix can be enhanced or not. We actually show
that for an R-matrix $R$ that satisfies our criterion there exists a unique $%
\mu $ such that $\left( PR,\mu \right) $ and $\left( RP,\mu ^{-1}\right) $
satisfy all the conditions $\left( \ref{YB3}\right) $, $\left( \ref{ENH1}%
\right) $, $\left( \ref{ENH3}\right) $ and $\left( \ref{ENH4}\right) $. As a
direct consequence our construction also gives an $\alpha $ and $\beta $
such that $\left( PR,\mu ,\alpha ,\beta \right) $ and $\left( RP,\mu
^{-1},\alpha ,\beta \right) $ are enhanced R-matrices in the sense of
definition \ref{e1}. We would like to stress that our $\mu ,\alpha $ and $%
\beta $ only depend on the R-matrix $R$, and not on any information coming
from Lie algebras. However, in order to put this problem in terms of
category theory, where we think it belongs, we have to assume that $R$ is 
{\it biinvertible}. This means that not only $R$ itself is invertible, but also
its second transpose $R^{t_2}$, where 
\[
\left( R^{t_2}\right) _{cd}^{ab}=R_{cb}^{ad}. 
\]

\noindent Although this seems to be a restriction, we will show in section 4
that any R-matrix that can be enhanced in the sense of def.\ref{e2} is
necessarily biinvertible. However, this is not true for R-matrices that can
be enhanced in the sense of def.\ref{e1}. A counterexample is the
permutation matrix $P$. It is easy to check that this matrix is not
biinvertible and, since $P^2=I$, we see that $\left( P^2,I,1,n\right) 
$, with $n^2$ the dimension of $P$, is an enhanced R-matrix in the sense of def.\ref{e1}. We don't know how
many other known not biinvertible R-matrices can be enhanced and if there is
a way to understand their enhancement. Despite this little gap in our
explanation of enhancement, we still think our results are interesting
enough. For example, in section 5 we show which $4\times 4$ R-matrices ($n=2$) can be enhanced and how.  
These R-matrices were classified by Hietarinto \cite{HIE} (see also \cite{MAJ}). It would be interesting to 
see what kind of knot invariants they give. It would also be very interesting to apply our results to
the R-matrices found by Cremmer and Gervais \cite{GER} and those classified
by Van den Hijligenberg \cite{NICO} and see if they can be enhanced and if
so, what kind of knot invariants they give. The latter tried to generalize
the work of Hazewinkel and classified a subclass of the R-matrices under the
restriction 
\[
R_{cd}^{ab}\neq 0\Longrightarrow \left\{ a,b\right\} =\left\{ c,d\right\} %
\smallskip\ \text{or}\smallskip\ a=\sigma \left( b\right) ,c=\sigma \left(
d\right) 
\]

\noindent with $\sigma \left( i\right) =n+1-i.$

This paper is organized as follows. In section 2 we recall the basic facts
about quasi-triangular Hopf algebras and braided categories. We followed 
\cite{KAS} closely and refer for the details and the proofs to the same. In section 3 we
explain the dual theory. This is the theory of dual quasi-triangular Hopf
algebras and braided categories. For the details and the proofs we refer to \cite
{KAS} and \cite{MAJ}. Section 4 contains our own results. In section 5 we
enhance all biinvertible $4\times 4$ R-matrices ($n=2$). Finally in the appendix we give
an elementary proof of our main theorem for the reader who is interested in
our results, but does not want to go through all the category theory in
sections 2, 3 and 4. 

\end{section}

\begin{section}{The basic idea}

\noindent Let $H$ be a Hopf algebra with comultiplication $\Delta$, counit $\varepsilon$ and invertible antipode $S.$ We say
that $H$ is quasi-triangular if there exists an invertible element ${\cal R}%
\in H\otimes H$ such that 
\begin{equation}  \label{RM1}
\Delta ^{{\rm op}}\left( x\right) ={\cal R}\Delta \left( x\right) {\cal R}^{-1},
\end{equation}
\begin{equation}  \label{RM2}
\left( \Delta \otimes {\rm id}\right) \left( {\cal R}\right) ={\cal R}_{13}%
{\cal R}_{23},
\end{equation}
\begin{equation}  \label{RM3}
\left( {\rm id}\otimes \Delta \right) \left( {\cal R}\right) ={\cal R}_{13}%
{\cal R}_{12}.
\end{equation}

\noindent Here $\Delta^{{\rm op}}$ is the opposite comultiplication defined by 
$\Delta^{{\rm op}}=\tau\Delta$, where $\tau$ is the flip $\tau(x\otimes y)=y\otimes x$. 
It is easy to check that these properties imply 
\begin{equation}
{\cal R}_{12}{\cal R}_{13}{\cal R}_{23}{\cal =R}_{23}{\cal R}_{13}{\cal R}%
_{12},  \label{RM4}
\end{equation}
\begin{equation}
\left( \varepsilon \otimes {\rm id}\right) \left( {\cal R}\right) =1=\left( 
{\rm id}\otimes \varepsilon \right) \left( {\cal R}\right) ,  \label{RM5}
\end{equation}
\begin{equation}
\left( S\otimes {\rm id}\right) \left( {\cal R}\right) ={\cal R}^{-1}=\left( 
{\rm id}\otimes S^{-1}\right) \left( {\cal R}\right) ,  \label{RM6}
\end{equation}
\begin{equation}
\left( S\otimes S\right) \left( {\cal R}\right) ={\cal R}.  \label{RM7}
\end{equation}
If such a ${\cal R}$ exists, then $\left( \ref{RM4}\right) $ shows that
every evaluation of ${\cal R}$ in a representation of $H$ satisfies the
QYB-equation $\left( \ref{YB}\right) $. That is why ${\cal R}$ is called the
universal R-matrix of $H.$ Now suppose $H$ is quasi-triangular with ${\cal R}%
=\sum r_i\otimes t_i.$ Then there is a well known lemma that says that $S^2$
is an inner automorphism of $H.$

\begin{lemma}\label{u}
Under the previous hypothesis, the elements $u=\sum
S\left( t_i\right) r_i$ and $v=S\left( u\right) $ are invertible elements in 
$H$ such that 
\begin{equation}
S^2\left( x\right) =uxu^{-1}=v^{-1}xv.  \label{ANTP1}
\end{equation}

\noindent The element $uv=vu$ is central in $H,$ and satisfies 
\begin{equation}
\Delta \left( uv\right) =\left( {\cal R}_{21}{\cal R}\right) ^{-2}\left(
uv\otimes uv\right)   \label{ANTP2}
\end{equation}
\end{lemma}

\noindent In order to see the connection with knot invariants we have to
consider knots as morphisms in a special tensor category ${\cal T}$, the
category of tangles. The objects of ${\cal T}$ are finite sequences of $+$
and $-$ signs. Their tensor product is the sequence that we get by putting
one sequence after the other. The identity object is just the empty
sequence. If we put two such sequences one above the other, then a
morphism in ${\cal T}$ is an equivalence class of oriented tangles between
them. By this we mean that, taking an arbitrary tangle in the equivalence
class, all $+$ and $-$ signs are either the head or the tail of a strand of this
tangle. The head of a strand is
attached to a $+$ sign if this head is pointing downward and it is attached
to a $-$ sign if the head is pointing upward. A tail of a strand is attached
to a $+$ sign if it is pointing downward and attached to a $-$ sign if it is
pointing upward. Two tangles are equivalent if one can be obtained from the
other by only applying homotopies of the tangle diagram and the Reidemeister
moves 1,2 and 3. In fig.\ref{f1}. we show an example.

\begin{figure}[htb]
\begin{center}
\mbox{\epsfxsize=0.3\hsize\epsfbox{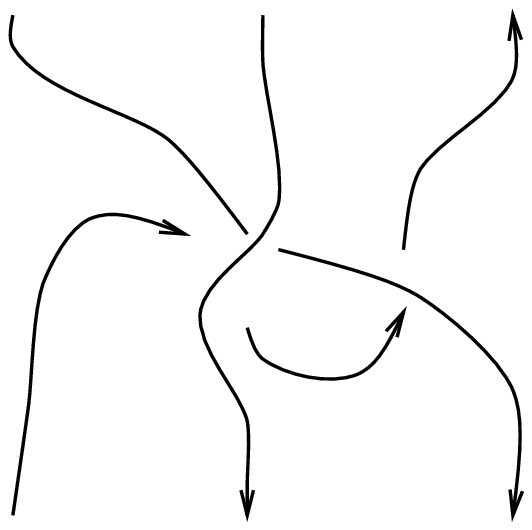}}
\caption{A tangle.}
\label{f1}
\end{center}
\end{figure}

\noindent The composition of two tangles is obtainded by putting the first
tangle on top of the second. This means that we 'read' the tangles from
bottom to top. The tensor product of two tangles is given by juxtaposition
and the identity endomorphism of an object is just the set of straight
vertical strands with orientation determined by the signs in the sequence.
It is an easy exercise to show that these objects and morphisms define a
tensor category, with the empty set being the identity object. In fact they
have more structure, which makes tangles intimately related to the
representation theory of quasi-triangular Hopf algebras. We can define
something called a {\it braiding} in the category of tangles.

\begin{definition}
A braiding in a tensor category $C$ is a family of natural isomorphisms $%
c_{V,W}:V\otimes W\rightarrow W\otimes V$, indexed by all pairs of objects $%
V,W\in {\rm Obj}\left( C\right) ,$ such that 
\begin{equation}
c_{U,V\otimes W}=\left( {\rm id}_V\otimes c_{U,W}\right) \left(
c_{U,V}\otimes {\rm id}_W\right)   \label{C1}
\end{equation}
\begin{equation}
c_{U\otimes V,W}=\left( c_{U,W}\otimes {\rm id}_V\right) \left( {\rm id}%
_U\otimes c_{V,W}\right)   \label{C2}
\end{equation}
\end{definition}

\noindent The naturality of $c$ means that 
\[
c\left( f\otimes g\right) =\left( g\otimes f\right) c, 
\]
where $f$ and $g$ are morphisms in $C$ and $c$ should be understood with the
right subscripts$.$ It is not difficult to see that $\left( \ref
{C1}\right) $ and $\left( \ref{C2}\right) $ imply the following identity 
\begin{equation}
\begin{array}{l}
\left( c_{V,W}\otimes {\rm id}_U\right) \left( {\rm id}_V\otimes
c_{U,W}\right) \left( c_{U,V}\otimes {\rm id}_W\right) = \\ 
\left( {\rm id}_W\otimes c_{U,V}\right) \left( c_{U,W}\otimes {\rm id}%
_V\right) \left( {\rm id}_U\otimes c_{V,W}\right) .
\end{array}
\label{C3}
\end{equation}

\noindent If we take $U=V=W,$ then identity \ref{C3} is similar to identity 
\ref{YB3}. That is why such a $c$ is sometimes called a Yang-baxter
operator on $V.$ The braiding in the category of tangles is defined in fig.\ref{f2}. 
The inverse of $c$ one gets by changing the overcrossings in undercrossings.

\begin{figure}[htb]
\begin{center}
\mbox{\epsfysize=6cm\epsfbox{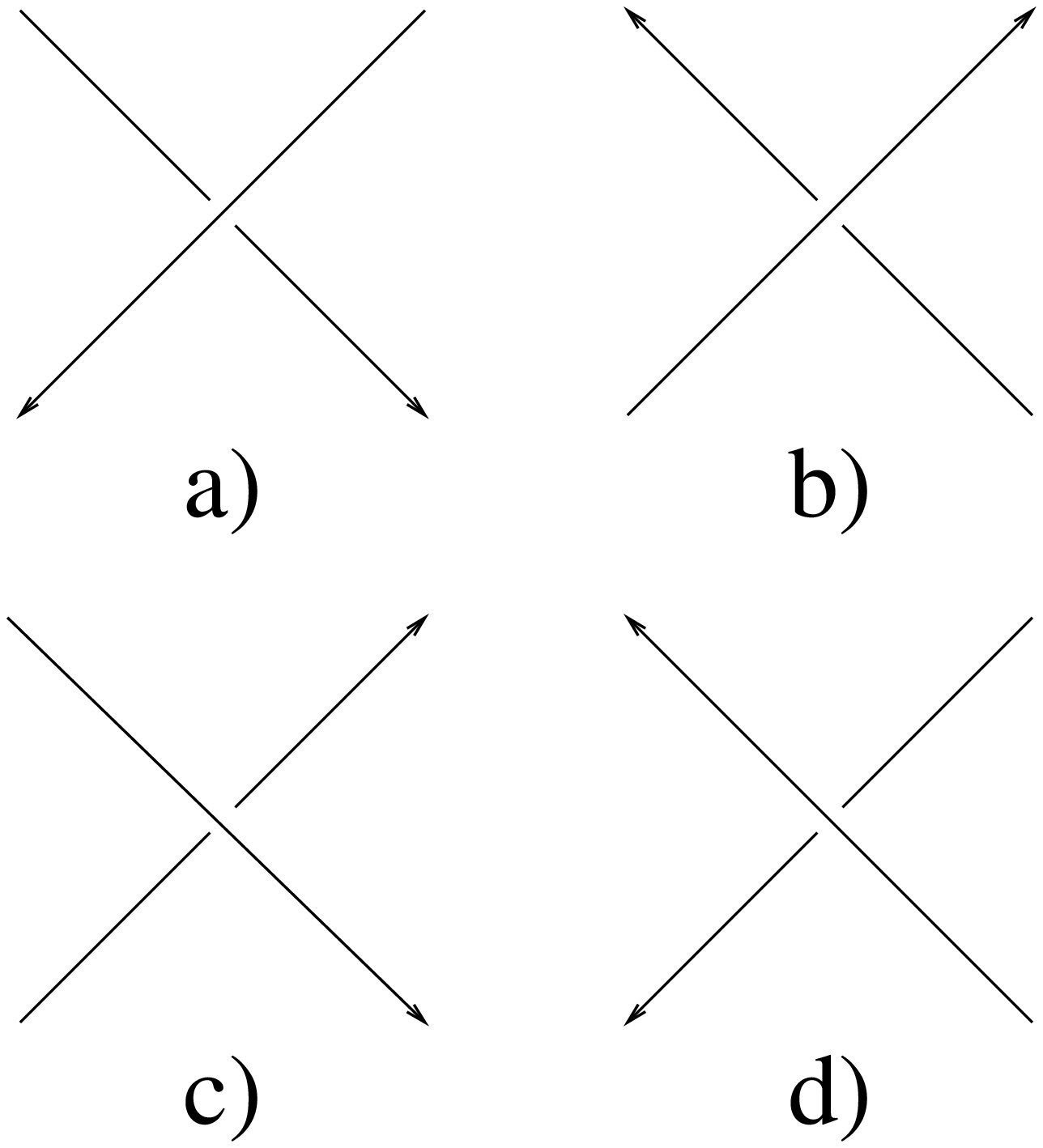}}
\caption{a) $c_{+,+}$ b) $c_{-,-}$ c) $c_{-,+}$ d) $c_{+,-}$.}
\label{f2}
\end{center}
\end{figure}

In the same way as certain algebras or groups can be presented by a finite set of
generators and relations between them there are categories that can be
presented by a finite set of generating morphisms and relations between
these morphisms. The category ${\cal T}$ can be presented
by the generators $X^{+},X^{-},\cup ,\cup ^{-},\cap ,\cap ^{-}$ (see fig.\ref
{f4}). For the defining relations see \cite{TUR3} or \cite{KAS}. 

\begin{figure}[hbt]
\begin{center}
\mbox{\epsfysize=6cm\epsfbox{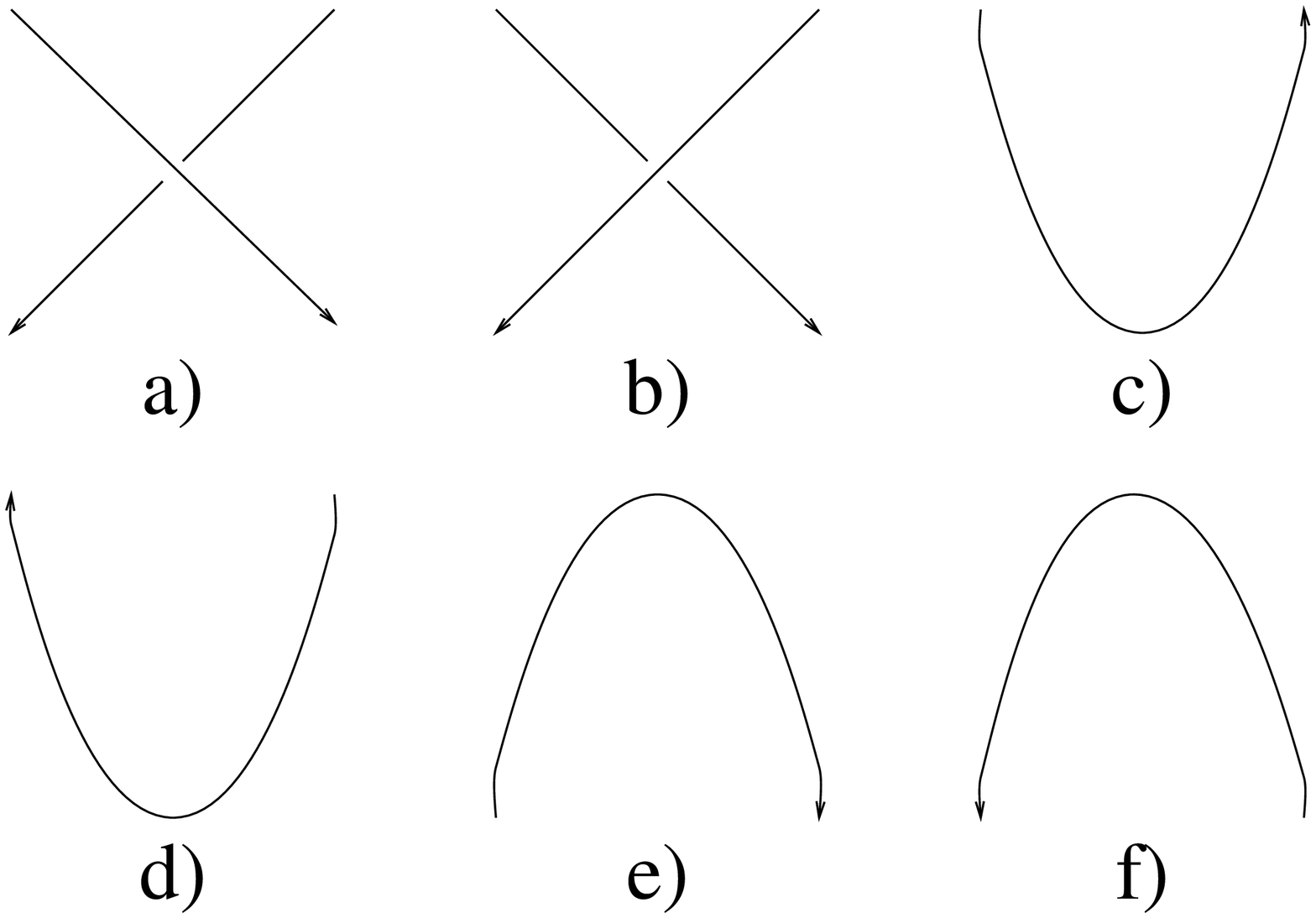}}
\caption{a) $X^+$  b) $X^-$ c) $\cup$ d) ${\cup}^-$ e) $\cap$ f) ${\cap}^-$.}
\label{f4}
\end{center}
\end{figure}

Another {\it braided} {\it category}
is the category of finite dimensional $H$-modules, where $H$ is a
quasi-triangular Hopf algebra (Kassel calls them {\it braided }Hopf
algebras, but the same name is used by Majid \cite{MAJ} for Hopf algebras in a
braided category, which is why we prefer to use the 'old' name). The theorem
is actually a bit stronger than that.

\begin{theorem}
\label{braid} Let $H$ be a Hopf algebra and $_H{\cal M}$ its
category of finite dimensional left modules. If $H$ is quasi-triangular, then%
{\it \ }$_H{\cal M}$ is a braided category. Conversely, if $_H{\cal M}$ is a
braided category and $H$ is finite dimensional, then $H$ is a
quasi-triangular Hopf algebra.
\end{theorem}

\noindent {\bf Sketch of a Proof. }Suppose $H$ is quasi-triangular and $%
{\cal R}$ is the universal R-matrix of $H$, then we define the braiding $c$
in $_H{\cal M}$ by 
\begin{equation}  \label{C4}
c_{V,W}\left( v\otimes w\right) =\tau _{V,W}\left( {\cal R}\left( v\otimes
w\right) \right) ,
\end{equation}

\noindent where $\tau _{V,W}$:$V\otimes W\rightarrow W\otimes V$ is the flip
operator. It is easy to check that this really defines a braiding.
Conversely, suppose that $_H{\cal M}$ is braided and $H$ is finite
dimensional. Then $H$ is a finite dimensional $H$-module with the action
defined by its multiplication. We define the invertible element 
\[
{\cal R}=\tau _{H,H}\left( c_{H,H}\left( 1\otimes 1\right) \right) . 
\]
\noindent It is easy to check that $H$ becomes quasi-triangular with
universal R-matrix ${\cal R}$\TeXButton{End Proof}{\hfill\endproof}

\medskip\ 

\noindent Instead of the category of left $H$-modules, we could also
consider the category of right $H$-modules ${\cal M}_H.$ The theorem above
then goes through as in the case of left modules, but, of course, we now
have to define 
\[
c_{V,W}\left( v\otimes w\right) =\tau _{V,W}\left( \left( v\otimes w\right) 
{\cal R}\right) . 
\]
If we consider framed tangles instead of normal tangles we get the category
of framed tangles ${\cal FT}$. Its objects are the same as those of ${\cal T%
},$ but the morphisms are equivalence classes of framed tangles. The framing
of a tangle is an equivalence class of normal vector fields on the tangle.
These framed tangles can be depicted as normal tangles with the convention
that the framing comes orthogonally out of the paper. With respect to the
orientation nothing changes but the equivalence relation has to be modified.
Two framed tangles are said to be equivalent if one can be obtained from the
other by isotopy of their diagrams, Reidemeister moves 2 and 3, and the
modified Reidemeister move 1' depicted in fig.\ref{f3}.

\begin{figure}[htb]
\begin{center}
\mbox{\epsfxsize=0.3\hsize\epsfbox{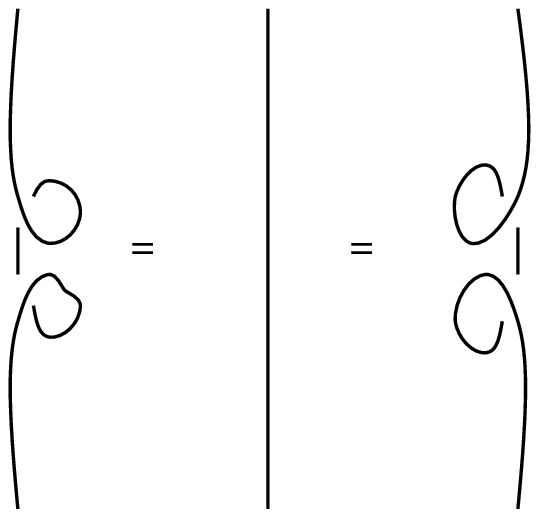}}
\caption{Reidemeister move 1'.}
\label{f3}
\end{center}
\end{figure}

\noindent This category is a tensor category in the same way as ${\cal T}$
and the braiding in ${\cal T}$ defines a braiding in ${\cal FT}$ as well. It
is also possible to define a last bit of extra structure on our category,
called {\it left duality }and {\it twist}.

\begin{definition}
A tensor category $C$ with identity object $I$ is a category with left
duality if for every object $V$ in $C$ there exists an object ${^{*}V}$ and
morphisms 
\[
b_V:I\rightarrow V\otimes {^{*}V}\smallskip\ {\rm and}\smallskip\
d_V:{^{*}V}\otimes V\rightarrow I
\]

\noindent in the category $C$ such that 
\begin{equation}
\left( {\rm id}_V\otimes d_V\right) \left( b_V\otimes {\rm id}_V\right) =%
{\rm id}_V\ \,{\rm and}\,\ \left( d_V\otimes {\rm id}_{^{*}V}\right) \left( 
{\rm id}_{^{*}V}\otimes b_V\right) ={\rm id}_{^{*}V}.  \label{dual1}
\end{equation}
\end{definition}

\noindent Using this definition we can define the {\it transpose }$%
{^{*}f}:{^{*}V}\rightarrow {^{*}U}$ of a morphism $f:U\rightarrow V$ in $C$ by 
\[
^{\ast }f=\left( d_V\otimes {\rm id}_{^{*}U}\right) \left( {\rm id}%
_{^{*}V}\otimes f\otimes {\rm id}_{^{*}U}\right) \left( {\rm id}%
_{^{*}V}\otimes b_U\right) . 
\]

\noindent In this way left duality defines a functor $*:C\rightarrow C.$
In general this functor is not an involution, which is why not every R-matrix 
gives a representation of the category of tangles. We will explain this in section 4.

\begin{theorem}
In a braided category with left duality there are natural
equivalences $u,v^{-1}\in {\rm Nat}\left( {\rm Id},*^2\right) $ defined by 
\[
u_V=\left( d_V\otimes {\rm id}\right) \left( c_{V,{^{*}V}}\otimes {\rm id}%
\right) \left( {\rm id}\otimes b_{^{*}V}\right) , 
\]
\[
u_V^{-1}=\left( {\rm id}\otimes d_{^{*}V}\right) \left( c_{{^{**}V},V}\otimes 
{\rm id}\right) \left( {\rm id}\otimes b_V\right) , 
\]
\[
v_V=\left( d_{^{*}V}\otimes {\rm id}\right) \left( {\rm id}\otimes
c_{V,{^{*}V}}\right) \left( {\rm id}\otimes b_V\right) , 
\]
\[
v_V^{-1}=\left( d_V\otimes {\rm id}\right) \left( {\rm id}\otimes
c_{{^{**}V},V}\right) \left( b_{^{*}V}\otimes {\rm id}\right) , 
\]
\noindent obeying 
\[
u_{V\otimes W}=c_{V,W}^{-1}c_{W,V}^{-1}\left( u_V\otimes u_W\right) , 
\]
\[
v_{V\otimes W}=c_{V,W}^{-1}c_{W,V}^{-1}\left( v_V\otimes v_W\right) , 
\]
and such that 
\[
{^{\ast *}}f=u_W\circ f\circ u_V^{-1}=v_W^{-1}\circ f\circ v_V,\,\,\,\,\forall
f:V\rightarrow W. 
\]
\end{theorem}

\noindent In general left duality is only unique up to an isomorphism. If $%
\left( \times ,b^{\times },d^{\times }\right) $ defines another left
duality, then one can define an isomorphism $f_V:^{\times }V\rightarrow
{^{*}V} $ for every object $V$ by 
\[
f_V=\left( d_V^{\times }\otimes {\rm id}\right) \left( {\rm id}\otimes
b_V\right) ,\,\,f_V^{-1}=\left( d_V\otimes {\rm id}\right) \left( {\rm id}%
\otimes b_V^{\times }\right) . 
\]

\noindent So we get 
\[
d_V^{\times }=d_V\left( f_V\otimes {\rm id}\right) ,\,\,b_V^{\times }=\left( 
{\rm id}\otimes f_{V\stackunder{}{\stackrel{}{}}}^{-1}\right) b_V. 
\]
\noindent There is a similar notion of right duality.

\begin{definition}
A tensor category $C$ with identity object $I$ is a category with right
duality if for every object $V$ in $C$ there exists an object $V^{*}$ and
morphisms 
\[
b_V^{^{\prime }}:I\rightarrow V^{*}\otimes V\smallskip\ {\rm and}\smallskip\
d_V^{^{\prime }}:V\otimes V^{*}\rightarrow I 
\]

\noindent in the category $C$ such that 
\begin{equation}
\left( {\rm id}_{V^{*}}\otimes d_V^{^{\prime }}\right) \left( b_V^{^{\prime
}}\otimes {\rm id}_{V^{*}}\right) ={\rm id}_{V^{*}}\ \,{\rm and}\,\ \left(
d_V^{^{\prime }}\otimes {\rm id}_V\right) \left( {\rm id}_V\otimes
b_V^{^{\prime }}\right) ={\rm id}_V.  \label{dual2}
\end{equation}
\end{definition}

\noindent Of course right duality is also unique up to isomorphism. If $H$
is a Hopf algebra with invertible antipode, then we can take the left dual $%
{^{*}V}={\rm Hom}(V,\Bbb C)$, the dual vector space, of every finite dimensional $H$%
-module $V$ and make it into a left $H$-module by defining 
\[
x\triangleright f\left( y\right) =f\left( S\left( x\right) y\right) 
\]

\noindent for every $f\in {^{*}V}$, every $x\in H$ and every $y\in V.$ The maps $d_V$ and $%
b_V$ are simply the evaluation and coevaluation map respectively. If $H$ is
quasi-triangular, we find the isomorphisms $u_V$ and $v_V$ to be the actions
of $u\in H$ and $v\in H$ on $V.$ 

The next notion we want to introduce is
that of a {\it twist.}

\begin{definition}
A twist in a braided tensor category $C$ with left duality is a family $%
\theta _V:V\rightarrow V$ of natural isomorphisms indexed by all objects $V$
in $C$ such that 
\begin{equation}
\theta _{V\otimes W}=\left( \theta _V\otimes \theta _W\right) c_{W,V}c_{V,W},
\label{t1}
\end{equation}
\begin{equation}
\theta _{^{*}V}={^{*}\left( \theta _V\right)}  \label{t2}
\end{equation}

\noindent for all objects $V,W$ in $C.$
\end{definition}

\noindent A braided tensor category with left duality and twist is called a 
{\it ribbon\ category}. In a ribbon category we also have right duality. We
just take $V^{*}={^{*}V}$ and define 
\begin{equation}
b_V^{^{\prime }}=\left( {\rm id}_{V^{*}}\otimes \theta _V\right)
c_{V,V^{*}}b_V,  \label{rd1}
\end{equation}
\begin{equation}
d_V^{^{\prime }}=d_Vc_{V,V^{*}}\left( \theta _V\otimes {\rm id}%
_{V^{*}}\right) .  \label{rd2}
\end{equation}

\noindent Before we give two examples of ribbon categories, we first state a
technical lemma that we will need in section 4.

\begin{lemma}
\label{twist} For any object $V$ in a ribbon category, we have 
\begin{equation}
\begin{array}{ll}
\theta _V^{-2} & =\left( d_V\otimes {\rm id}_V\right) \left( {\rm id}%
_{^{*}V}\otimes c_{V,V}^{-1}\right) \left( c_{V,{^{*}V}}b_V\otimes {\rm id}%
_V\right)  \\ 
& =\left( d_Vc_{V,{^{*}V}}\otimes {\rm id}_V\right) \left( {\rm id}_V\otimes
c_{V,{^{*}V}}b_V\right)  \\ 
& =\left( {\rm id}_V\otimes d_Vc_{V,{^{*}V}}\right) \left(
c_{V,{^{*}V}}^{-1}\otimes {\rm id}_{^{*}V}\right) \left( {\rm id}_V\otimes
b_V\right) 
\end{array}
\label{Theta}
\end{equation}
\end{lemma}

\noindent Of course our examples of ribbon categories are again the
category of framed tangles and the category of $H$-modules. We first
define the duality and the twist in ${\em {\cal FT}}.$ (Framed tangles are
often called ribbons, so that is where the name {\it ribbon category} comes
from.) Let $\varepsilon =\left( \varepsilon _1,\ldots ,\varepsilon _n\right) 
$ be a finite sequence of $+$ and $-$ signs. Define the dual object $%
{^{*}\varepsilon} =\left( -\varepsilon _n,\ldots ,-\varepsilon _1\right) .$
The morphisms $b_\varepsilon :\emptyset \rightarrow \varepsilon \otimes
{^{*}\varepsilon} $ and $d_\varepsilon :{^{*}\varepsilon} \otimes \varepsilon
\rightarrow \emptyset $ are the framed tangles depicted in fig.\ref{f5}
, their orientation being
completely determined by the signs in $\varepsilon .$ The relations $\left( 
\ref{dual1}\right) $ are easy to check in this case. Note that the transpose 
${^{*}L}$ of a tangle $L$ is obtained by rotation of the whole diagram through
an angle $\pi .$ 

\begin{figure}[htb]
\begin{center}
\mbox{\epsfysize=3cm\epsfbox{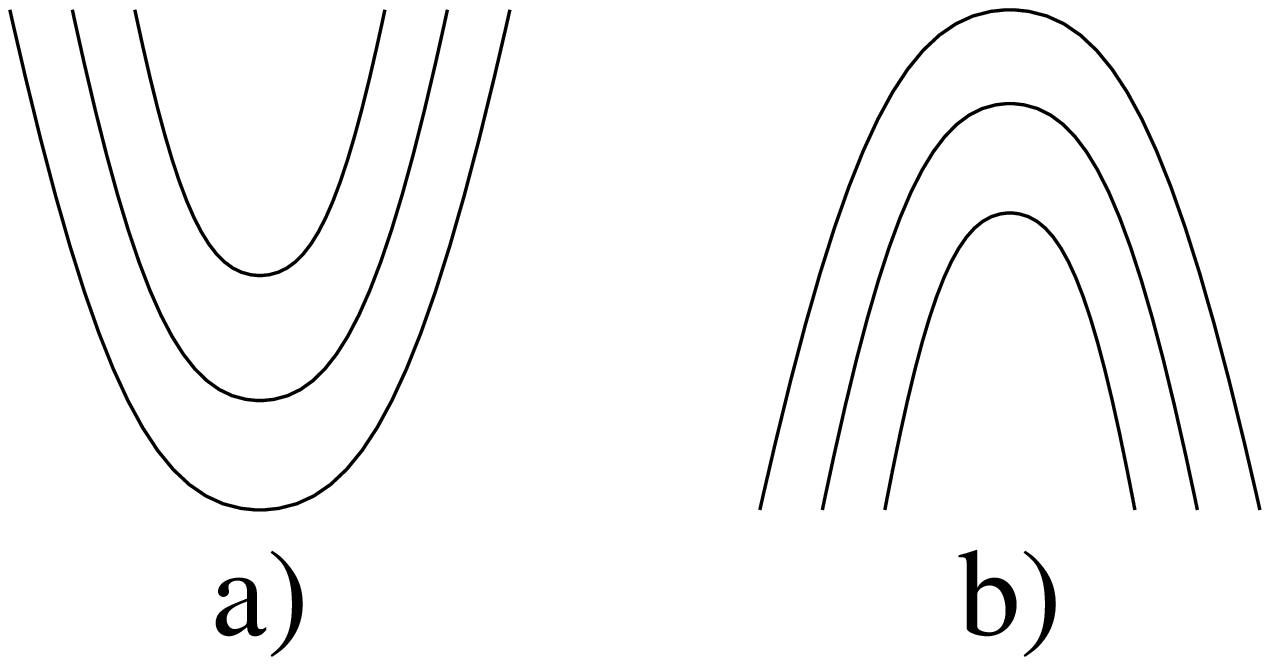}}
\caption{a) $b_{\epsilon}$ b) $d_{\epsilon}$.}
\label{f5}
\end{center}
\end{figure}

\noindent The twist on $+,$ denoted by $\varphi ,$ we define in fig.\ref{f5}.
The twist on an arbitrary object is then defined by the formulas $\left( \ref{t1}\right) $
and $\left( \ref{t2}\right) $. The category ${\em {\cal FT}}$ can be
presented by $X^{+},X^{-},\cup ,\cup ^{-},\cap ,\cap ^{-}$ and all the same
defining relations as in ${\cal T}$ except one: the one that corresponds to
Reidemeister move 1. This one has to be replaced by a relation that
expresses Reidemeister move 1'. The exact definition of this relation in $%
{\em {\cal FT}}$ we leave to the reader. 

\begin{figure}[htb]
\begin{center}
\mbox{\epsfysize=4cm\epsfbox{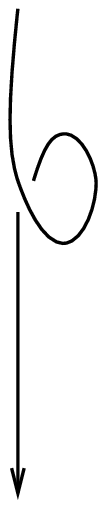}}
\caption{${\varphi}_+$.}
\label{f6}
\end{center}
\end{figure}

\noindent Now we can formulate a very
important property of ${\em {\cal FT}},$ called the {\it universality
property}

\begin{theorem}
\label{br} Let $C$ be a ribbon category and $V$ an object in $C.$
Then there exists a unique functor $F_V:{\cal FT}\rightarrow C$ preserving
the braiding, the left duality and the twist, such that $F_V\left( +\right)
=V$ and $F_V\left( -\right) ={^{*}V}.$
\end{theorem}

\noindent The functor $F_V$ has the following properties: 
\[
F\left( X^{+}\right) =c_{V,V},\,\,F\left( \varphi \right) =\theta
_V,\,\,F\left( \cup \right) =b_V,\,\,F\left( \cap \right) =d_V; 
\]
\[
F\left( X^{-}\right) =c_{V,V}^{-1},\,\,F\left( T^{+}\right)
=c_{V,{^{*}V}}^{-1},\,\,F\left( T^{-}\right) =c_{{^{*}V},V}, 
\]
\[
F\left( Y^{+}\right) =c_{{^{*}V},V}^{-1},\,\,F\left( Y^{-}\right)
=c_{V,{^{*}V}}, 
\]
\[
F\left( Z^{+}\right) =c_{{^{*}V},{^{*}V}},\,\,F\left( Z^{-}\right)
=c_{{^{*}V},{^{*}V}}^{-1},\,\,F\left( \varphi ^{-}\right) =\theta _V^{-1}. 
\]

\noindent The tangles $Z^{+},Y^{+},T^{+}$ are depicted in fig.\ref{f6}. Their inverses one gets by changing overcrossings in undercrossings. 

\begin{figure}[htb]
\begin{center}
\mbox{\epsfysize=3cm\epsfbox{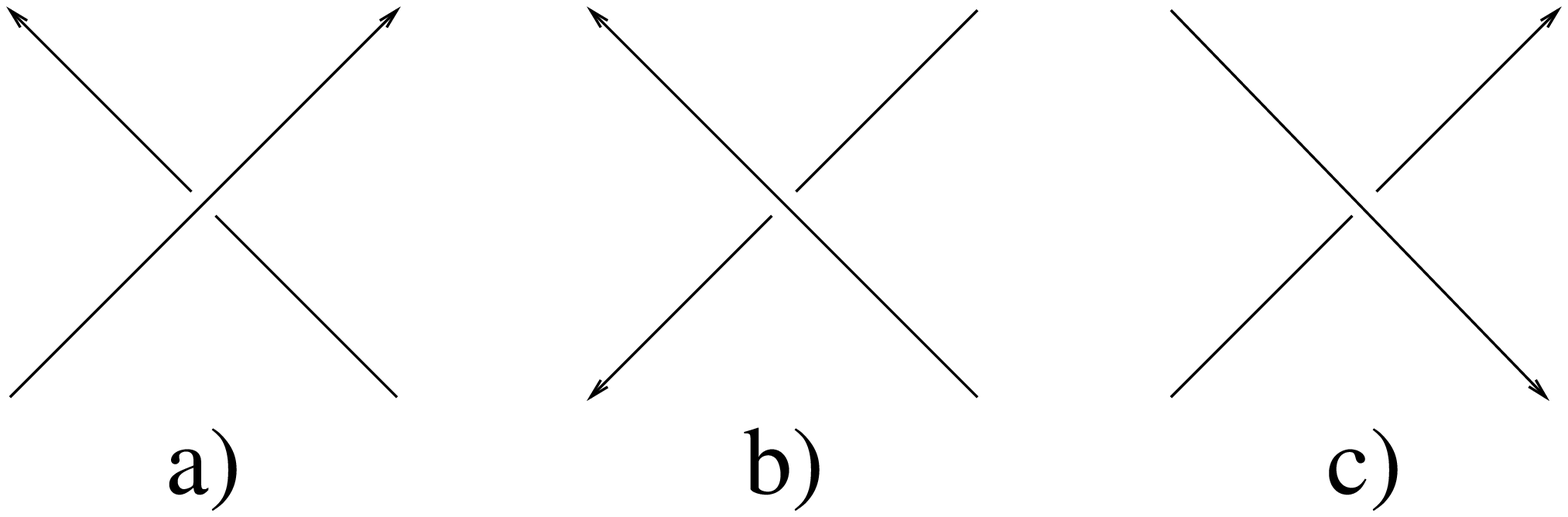}}
\caption{a) $Z^+$ b) $Y^+$ c) $T^+$.}
\label{f7}
\end{center}
\end{figure}

\noindent The values
on $\cup ^{-}$ and $\cap ^{-}$can be easily computed from the formulas 
\[
\cup ^{-}=\left( \uparrow \otimes \varphi ^{-}\right) \circ Y^{+}\circ \cup
_{,} 
\]
\[
\cap ^{-}=\cap \circ Y^{-}\circ \left( \varphi \otimes \uparrow \right) . 
\]

\noindent From $\left( \ref{rd1}\right) $ and $\left( \ref{rd2}\right) $ we
get $F\left( \cup ^{-}\right) =b_V^{^{\prime }}$ and $F\left( \cap
^{-}\right) =d_V^{^{\prime }}.$ It follows that $X^{\pm },Y^{\pm },\varphi
^{\pm },\cup ,\cap $ is another set of generators for ${\cal FT}$. Conversely, we can
also express $Z^{\pm },Y^{\pm },T^{\pm }$ in terms of $X^{\pm },\cup ^{\pm
},\cap ^{\pm }.$

\begin{lemma}
\label{rel} The following relations hold in the categories ${\cal T%
}$ and ${\cal FT}$: 
\[
Y^{\pm }=\left( \uparrow \downarrow \cap ^{-}\right) \left( \uparrow X^{\pm
}\uparrow \right) \left( \cup ^{-}\downarrow \uparrow \right) ,
\]
\[
T^{\pm }=\left( \cap \downarrow \uparrow \right) \left( \uparrow X^{\pm
}\uparrow \right) \left( \uparrow \downarrow \cup \right) ,
\]
\[
Z^{\pm }=\left( \cap \uparrow \uparrow \right) \left( \uparrow \cap
\downarrow \uparrow \uparrow \right) \left( \uparrow \uparrow X^{\pm
}\uparrow \uparrow \right) \left( \uparrow \uparrow \downarrow \cup \uparrow
\right) \left( \uparrow \uparrow \cup \right) ,
\]
\[
Z^{\pm }=\left( \uparrow \uparrow \cap ^{-}\right) \left( \uparrow \uparrow
\downarrow \cap ^{-}\uparrow \right) \left( \uparrow \uparrow X^{\pm
}\uparrow \uparrow \right) \left( \uparrow \cup ^{-}\downarrow \uparrow
\uparrow \right) \left( \cup ^{-}\uparrow \uparrow \right) .
\]
\end{lemma}

\noindent Now suppose that we have a ribbon category $C$ and a specified
object $V$ in $C$. Then we get an invariant of framed links with values in
End$\left( I\right) $, with $I$ the identity object in $C,$ if we consider
framed links as morphisms in ${\em {\cal FT}}$ from the empty set to itself.
If $C$ is a subcategory of the category ${\cal V_{{\Bbb C}}}$, the category
of complex vector spaces, we get an invariant with values in ${\bf {\Bbb C}}$%
. An important question is now whether or not we can find such ribbon categories,
 because all the machinery above only gives us a concrete framed knot
invariant if we come up with a concrete ribbon category different from ${\em 
{\cal FT}}$. Before we show an example we first give one more definition. 

Let $%
H$ be a braided Hopf algebra with invertible antipode $S.$ We defined in
lemma \ref{u} the central elements $u,v$ such that $S^2\left( x\right)
=uxu^{-1}=v^{-1}xv$.

\begin{definition}
A quasi-triangular Hopf algebra $H$ is called a ribbon algebra if there
exists an invertible central element $\theta \in H,$ that we call the ribbon%
{\it \ }element, such that 
\begin{equation}
\theta ^2=vu,\,\,\Delta \left( \theta \right) =\left( {\cal R}_{21}{\cal R}%
\right) ^{-1}\left( \theta \otimes \theta \right) ,\,\,\varepsilon \left(
\theta \right) =1,\,\,S\left( \theta \right) =\theta .  \label{rib1}
\end{equation}
\end{definition}

\noindent This definition is exactly the one needed to make the following
theorem hold.
\begin{theorem}
\label{r}For any ribbon algebra $H$, the tensor category $_H{\cal %
M}$ is a ribbon category with twist $\theta _V$ given on any $H$-module $V$
by the multiplication by the inverse of the ribbon element. \\ Conversely,
if $H$ is a finite-dimensional braided Hopf algebra and the braided category 
$_H{\cal M}$ is a ribbon category, then $H$ is a ribbon algebra with ribbon
element defined by $\theta =\theta _H\left( 1\right) ^{-1}.$
\end{theorem}

\noindent Of course there is a similar theorem for ${\cal M}_H$. Note that a
priori there may be more than one ribbon element in the same
quasi-triangular Hopf algebra.

To give a first example of a ribbon algebra we consider Sweedler's four
dimensional Hopf algebra ${\cal S}.$ This Hopf algebra is generated by 
\[
1,\,\,x,\,\,y 
\]
\noindent with defining relations 
\[
x^2=1,\,\,y^2=0,\,\,xy+yx=0. 
\]
\noindent The Hopf algebra structure is given by 
\[
\Delta \left( x\right) =x\otimes x,\,\,\Delta \left( y\right) =1\otimes
y+y\otimes 1,\,\,\varepsilon \left( x\right) =1,\,\,\varepsilon \left(
y\right) =0 
\]
\noindent and 
\[
S\left( x\right) =x,\,\,S\left( y\right) =xy. 
\]
\noindent It is easy to verify that for every $\lambda \in {\bf {\Bbb C}}$
the following expression defines an universal R-matrix of ${\cal S}$%
\begin{eqnarray*}
{\cal R}_\lambda &=&\frac 12\left( 1\otimes 1+1\otimes x+x\otimes 1-x\otimes
x\right) + \\
&&\frac \lambda 2\left( y\otimes y+y\otimes xy+xy\otimes xy-xy\otimes
y\right) .
\end{eqnarray*}
\noindent It is also easy to verify that ${\cal R}_\lambda ^{-1}=\tau _{%
{\cal S},{\cal S}}\left( {\cal R}_\lambda \right) .$ After a simple
calculation we get $u=v=x.$ This shows that $\theta =1$ defines a ribbon
element in ${\cal S}$. (Majid \cite{MAJ} Example 2.1.11 takes $x$ as a
ribbon element, but $x$ is not central.) So for every $\lambda $ we get a
ribbon category with trivial twist.

Using theorem \ref{r} we can show where to find the ribbon categories that
we need in order to construct the so called quantum invariants. Let {\em g}
be a semi-simple complex Lie algebra of finite dimension and $U({\em g})$
its universal enveloping algebra. Drinfeld \cite{DRI85, DRI87} defined a quantization of this algebra by introducing a formal parameter $h$
in the defining relations between the generators of $U({\em g})$ that
destroys its commutativity and cocommutativity. We denote this 'standard'
quantization, which is defined over the ring ${\bf {\Bbb C}}\left[ \left[
h\right] \right] ,$ by $U_h({\em g}).$ Drinfeld also proved that $U_h({\em g})$ is
a quasi-triangular Hopf algebra. Turaev and Reshetikhin \cite{RT} showed
that $U_h({\em g})$ is actually a ribbon algebra by proving that the square
root of $vu\in U_h({\em g})$ exists and that it satisfies the properties $%
\left( \ref{rib1}\right) $. Theorem \ref{r} shows that $_{U_h({\em g})}{\cal %
M}$ is a ribbon category, so by the universality property there exists for
every module $V$ in $_{U_h({\em g})}{\cal M}$ a unique functor $F_V:{\cal FT}%
\rightarrow _{U_h({\em g})}{\cal M}$ such that $F_V\left( +\right) =V$ and $%
F_V\left( -\right) ={^{*}V}.$ As the functor is defined on equivalence classes
of framed tangles, we get an invariant of framed links when we restrict our
functor to this particular kind of framed tangles. It can be shown that, if $%
V$ is a finite dimensional {\em g}-module, there exists a unique $U_h({\em g}%
)$-module $\tilde V$ of finite rank, defined over ${\bf {\Bbb C}}\left[
\left[ h\right] \right] ,$ such that $\tilde V\equiv V{\rm mod\,}h.$ So,
given a semi-simple finite dimensional complex Lie algebra {\em g} and a
finite dimensional {\em g}-module $V,$ we can define a unique framed knot
and link invariant that essentially comes from the quantization of $U({\em g}%
).$ That is why these invariants are called quantum invariants. Of course we
could choose a right module $W,$ instead of a left module. Then there is a
unique functor $F_W:{\cal FT}\rightarrow {\cal M}_{U_h({\em g})}$ such that $%
F_W\left( +\right) =W$ and $F_W\left( -\right) ={^{*}W}.$

There are now two questions that we want to discuss. The first is how to
derive invariants of ordinary knots and links from this beautifully designed
theory. The second is how do enhanced R-matrices fit into this picture. This
last question is especially interesting because not all known R-matrices
were found by evaluating the universal R-matrix of a quasi-triangular Hopf
algebra (see for example \cite{HAZ}, \cite{NICO}).

The first problem can be solved easily. If you look closely at all
arguments, you see that the whole theory would work for ordinary tangles as
well, if only we required the twist on $+$ in ${\cal T}$ to be trivial and
the twist on a specified object in another ribbon category to be trivial as
well.

\begin{lemma}
The category of tangles ${\cal T}$ is a ribbon category if we define the
twist in the following way: 
\[
\theta _{+}\left( \downarrow \right) =\downarrow \,\,{\rm and}%
\,\,\theta _{-}\left( \uparrow \right) =\uparrow .
\]

\noindent The definition of the twist is extended to all other tangles by
applying formulas $\left( \ref{t1}\right) $ and $\left( \ref{t2}\right) $.
\end{lemma}

\noindent {\bf Proof}. Trivial.

\begin{theorem}
\label{univ2} Let $C$ be a ribbon category and $V$ an object in $C$ such
that $\theta _V ={\rm id}_V.$ Then there exists a unique
functor $F_V:{\em {\cal T}}\rightarrow C$, preserving braiding, duality and
twist, that sends $+$ to $V$ and $-$ to ${^{*}V}.$
\end{theorem}

\noindent {\bf Proof}. The proof of this theorem is identical to that of
theorem \ref{br}, apart from the details concerning the twist. Since we
require $F_V$ to preserve the twist, the twist in ${\em {\cal T}}$ defined
in the lemma above, is exactly the one that makes the construction of $F_V$
possible.\TeXButton{End Proof}{\hfill\endproof}

\medskip\ 

\noindent So we have to know where to find this particular kind of ribbon
categories. If $H$ is a ribbon algebra with trivial ribbon element, then of
course its category of representations $_H{\cal M}$ is such a category.
Thus the condition $vu=1$ seems a good one, if we already know that $H$ is a
quasi-triangular Hopf algebra. However there is one problem with this
condition. It is not a condition on R-matrices but rather on
quasi-triangular Hopf algebras, so enhancement does not seem to come into
the picture. Nevertheless the next lemma shows that we are on the right way,
by making a link to what one might call the 'dual' theory. This lemma follows 
more or less directly from \cite{MAJ} Prop. 4.2.2. Nevertheless, we have included 
a proof, because we could not find the lemma in the literature. 

\begin{lemma}
\label{biin}Let $H$ be a quasi-triangular Hopf algebra with universal
R-matrix ${\cal R}${\em .} Let $M$ be a finite dimensional $H$-module, so
that $\rho _M\left( {\cal R}\right) =R$ is the R-matrix on $M\otimes M.$ This
R-matrix is biinvertible (see introduction). If we define 
\[
\tilde R=\left( \left( R^{t_2}\right) ^{-1}\right) ^{t_2}, 
\]

\noindent then we get $U_j^i=\rho _M\left( u\right) _j^i=\tilde R_{ja}^{ai}$
and $V_j^i=\rho _M\left( v\right) _j^i=\tilde R_{aj}^{ia}.$
\end{lemma}

\noindent {\bf Proof}. We prove that $R$ is biinvertible by showing that 
\[
\rho _M\left( \left( {\rm id}\otimes S\right) {\cal R}\right) =\tilde R. 
\]

\noindent Let ${\cal R}=\sum s_i\otimes t_i$. Then we get 
\[
\begin{array}{cl}
1 & =\left( {\rm id}\otimes \varepsilon \right) \left( {\cal R}\right)
=\left( {\rm id}\otimes \left( S*{\rm id}\right) \right) \left( {\cal R}%
\right) =\left( {\rm id}\otimes \mu \right) \left( {\rm id}\otimes S\otimes 
{\rm id}\right) \left( {\rm id}\otimes \Delta \right) \left( {\cal R}\right)
\\ 
& =\left( {\rm id}\otimes \mu \right) \left( {\rm id}\otimes S\otimes {\rm id%
}\right) \left( {\cal R}_{13}{\cal R}_{12}\right) =\sum s_is_j\otimes
S\left( t_j\right) t_i.
\end{array}
\]

\noindent So we see 
\[
\begin{array}{cl}
\delta _c^a\delta _d^b & =\rho _M\left( \sum s_is_j\otimes S\left(
t_j\right) t_i\right) _{cd}^{ab}=\sum \rho _M\left( s_i\right) _\alpha
^a\rho _M\left( s_j\right) _c^\alpha \rho _M\left( S\left( t_j\right)
\right) _\beta ^b\rho _M\left( t_i\right) _d^\beta \\ 
& =\rho _M\left( \sum s_i\otimes t_i\right) _{\alpha d}^{a\beta }\rho
_M\left( \sum s_j\otimes S\left( t_j\right) \right) _{c\beta }^{\alpha
b}=R_{\alpha d}^{a\beta }\tilde R_{c\beta }^{\alpha b}.
\end{array}
\]

\noindent The proof of the other two assertions now follows easily: 
\[
\begin{array}{l}
U_b^a=\rho _M\left( \sum S\left( t_i\right) s_i\right) _b^a=\tilde
R_{b\alpha }^{\alpha a}, \\ 
\\ 
V_b^a=\rho _M\left( \sum s_iS\left( t_i\right) \right) _b^a=\tilde R_{\alpha
b}^{a\alpha }.\TeXButton{End Proof}{\hfill\endproof}
\end{array}
\]

\noindent The key observation here is that the matrices $U$ and $V$ can be
derived directly from the R-matrix $R$, without using any information coming
from $H.$ So, if one is optimistic, one hopes that the condition $VU=I,$ for
any biinvertible R-matrix $R,$ is the one equivalent to enhancement.

\end{section}

\begin{section}{The dual theory}

\noindent 
In the previous section we wrote that lemma \ref{biin} defined a link with
the 'dual' theory. In this section we explain what we mean by that.

\begin{definition}
Let $R$ be an $n^2\times n^2$ R-matrix. If $K\left\langle T\right\rangle =K\left\langle
T_1^1,\ldots ,T_n^n\right\rangle $ is the free bialgebra on $n^2$
generators, defined by 
\[
\Delta \left( T_j^i\right) =T_a^i\otimes T_j^a,\smallskip\ \varepsilon
\left( T_j^i\right) =\delta _j^i, 
\]

\noindent then we define the bialgebra $A(R)$ to be 
\[
K\left\langle T\right\rangle /\left( RT_1T_2-T_2T_1R\right) . 
\]
\end{definition}

\noindent For the proof that $J=\left( RT_1T_2-T_2T_1R\right) $ is really a
bialgebra ideal see for example \cite{HAZ2}. This bialgebra is not
quasi-triangular, but it has a 'dual' property.

\begin{definition}
Let $H$ be a bialgebra or a Hopf algebra. We say that $H$ is dual
quasi-triangular if there exists a convolution-invertible map ${\em R}%
:H\otimes H\rightarrow {\bf {\Bbb C}}$ such that 
\begin{equation}
\sum b_{\left( 1\right) }a_{\left( 1\right) }{\em R}\left( a_{\left(
2\right) }\otimes b_{\left( 2\right) }\right) =\sum {\em R}\left( a_{\left(
1\right) }\otimes b_{\left( 1\right) }\right) a_{\left( 2\right) }b_{\left(
2\right) },  \label{dr1}
\end{equation}
\begin{equation}
\begin{array}{c}
{\em R}\left( ab\otimes c\right) =\dsum {\em R}\left( a\otimes c_{\left(
1\right) }\right) {\em R}\left( b\otimes c_{\left( 2\right) }\right) , \\ 
{\em R}\left( a\otimes bc\right) =\dsum {\em R}\left( a_{\left( 1\right)
}\otimes c\right) {\em R}\left( a_{\left( 2\right) }\otimes b\right) .
\end{array}
\label{dr2}
\end{equation}

\noindent for all $a,b,c\in H.$
\end{definition}

\noindent Here we used Sweedler's notation: $\Delta \left( x\right) =\sum
x_{\left( 1\right) }\otimes x_{\left( 2\right) }.$ One can check that $A(R)$
is a dual quasi-triangular bialgebra with 
\[
{\em R}\left( T_c^a\otimes T_d^b\right) =R_{cd}^{ab}. 
\]

\noindent In a way $A(R)$ is the universal bialgebra with dual
quasi-triangular structure defined by $R.$

\begin{theorem}
Let $H$ be a dual quasi-triangular bialgebra with $n^2$
generators $\left\{ x_1^1,\ldots ,x_n^n\right\} $ and 1, having the
coalgebra structure given by 
\[
\Delta \left( x_j^i\right) =x_a^i\otimes x_j^a,\smallskip\ \varepsilon
\left( x_j^i\right) =\delta _{j.}^i 
\]
\noindent If 
\[
{\em R}\left( x_a^i\otimes x_b^j\right) =R_{ab}^{ij}, 
\]

\noindent then $R$ is an R-matrix and $H$ is a quotient of $A(R).$
\end{theorem}

\noindent In general $A(R)$ does not define a Hopf algebra. Even after
dividing out by some extra relations we do not always get a Hopf algebra.
However, if $R$ is biinvertible, then Majid \cite{MAJ} showed that there is
a formal way to extend $A(R)$ to a Hopf algebra $H(R).$

\begin{theorem}
a) Suppose that we can add relations to $A(R)$ such that the dual
quasitriangular structure descends to the quotient and gives us a dual
quasi-triangular Hopf algebra. Then 
\[
{\em R}\left( T\otimes T^{-1}\right) =\tilde R 
\]

\noindent obeys 
\[
\tilde R_{aj}^{ib}R_{lb}^{ak}=\delta _l^i\delta _j^k=R_{aj}^{ib}\tilde
R_{lb}^{ak},\smallskip\ \text{i.e. }\tilde R=\left( \left( R^{t_2}\right)
^{-1}\right) ^{t_2}. 
\]

\noindent Moreover, 
\[
S^2T=V^{-1}TV=UTU^{-1};\smallskip\ V_j^i=\tilde R_{aj}^{ia},\,\,U_j^i=\tilde
R_{ja}^{ai}, 
\]
\[
V_2=RV_2\tilde R,\smallskip\ V_1=\tilde RV_1R. 
\]

\noindent b) If $R$ is biinvertible, then we can enlarge $A(R)$ to obtain a
Hopf algebra $H(R),$ with the same dual quasi-triangular structure, by
adding formally the generators $T^{-1}=\left( \left( T^{-1}\right)
_j^i\right) _{i,j=1}^n$, with coalgebra structure 
\[
\Delta \left( T^{-1}\right) =\left( T^{-1}\otimes T^{-1}\right) ,\smallskip\
\varepsilon \left( T^{-1}\right) =I, 
\]

\noindent and additional relations 
\[
TT^{-1}=T^{-1}T=I, 
\]
\[
RT_1=T_2T_1RT_2^{-1},\smallskip\ T_2^{-1}RT_1T_2=T_1R,\smallskip\
T_1^{-1}T_2^{-1}R=RT_2^{-1}T_1^{-1}. 
\]

\noindent The antipode we define by 
\[
ST=T^{-1},\smallskip\ ST^{-1}=V^{-1}TV=UTU^{-1}, 
\]

\noindent and the dual quasi-triangular structure by 
\[
{\em R}\left( T\otimes T\right) ={\em R}\left( T^{-1}\otimes T^{-1}\right)
=R,\,\,{\em R}\left( T^{-1}\otimes T\right) =R^{-1},\,\,{\em R}\left(
T\otimes T^{-1}\right) =\tilde R. 
\]
\end{theorem}

\noindent The proof of part a) in the previous theorem can be found in \cite
{MAJ} Prop. 4.2.2. Part b) follows more or less automatically from Majid's
comments after this proposition. The universality of $A(R)$ and the theorem
above show us the following.

\begin{corollary}
\label{un}Let $H$ be a dual quasi-triangular Hopf algebra with $2n^2$
generators $\left\{ x_1^1,\ldots ,x_n^n,y_1^1,\ldots ,y_n^n\right\} $ and 1,
satisfying the relations 
\[
XY=I,\smallskip\ X_j^i=x_j^i,\,\,Y_j^i=y_j^i. 
\]

\noindent Suppose that $H$ has the coalgebra structure given by 
\[
\Delta \left( X\right) =X\otimes X,\,\,\Delta \left( Y\right) =Y\otimes
Y,\,\,\varepsilon \left( X\right) =\varepsilon \left( Y\right) =I. 
\]

\noindent If 
\[
{\em R}\left( X\otimes X\right) =R, 
\]

\noindent then $R$ is biinvertible and $H$ is a quotient of $H(R).$
Furthermore 
\[
{\em R}\left( Y\otimes X\right) =R^{-1},\,\,{\em R}\left( X\otimes Y\right)
=\tilde R,\,\,{\em R}\left( Y\otimes Y\right) =R, 
\]
\[
S\left( X\right) =Y,\,\,S\left( Y\right) =V^{-1}XV=UXU^{-1}. 
\]
\end{corollary}

\noindent In this context theorem \ref{braid} becomes

\begin{theorem}
Let $H$ be a Hopf algebra. If $H$ is dual quasi-triangular, then
the category of finite dimensional left $H$-comodules $^H{\cal M}$ is a
braided category. Conversely, if $H$ is finite dimensional and $^H{\cal M}$
is a braided category, then $H$ is a dual quasi-triangular Hopf algebra.
\end{theorem}

\noindent {\bf Sketch of a proof}. If $H$ is dual quasi-triangular, then $^H%
{\cal M}$ becomes a braided category with braiding 
\[
c_{V,W}\left( v\otimes w\right) =\sum {\em R}\left( v^{\left( 1\right)
}\otimes w^{\left( 1\right) }\right) w^{\left( 2\right) }\otimes v^{\left(
2\right) }, 
\]

\noindent for all $H$-modules $V$ and $W,$ where we define the comodule
structure explicitly by $\Delta _V\left( v\right) =\sum v^{\left( 1\right)
}\otimes v^{\left( 2\right) }.$ If $H$ is finite dimensional and $^H{\cal M}$
is a braided category, then we define the dual quasi-triangular structure on 
$H$ by 
\[
{\em R}\left( g\otimes h\right) =\left( \varepsilon \otimes \varepsilon
\right) \tau _{H,H}c_{H,H}\left( g\otimes h\right) .\TeXButton{End Proof}
{\hfill\endproof} 
\]

\medskip\ 

\noindent Of course there is a similar theorem for the category of right
comodules ${\cal M}^H.$ In that case the braiding is defined by 
\[
c_{V,W}\left( v\otimes w\right) =\sum w^{\left( 1\right) }\otimes v^{\left(
1\right) }{\em R}\left( v^{\left( 2\right) }\otimes w^{\left( 2\right)
}\right) . 
\]

\noindent We see that, if $R$ is biinvertible, the category $^{H\left(
R\right) }{\cal M}$ is a braided category. If $V$ is the left fundamental
corepresentation $V,$ with coaction $e^j\mapsto T_a^j\otimes e^a,$ then the
braiding is defined by 
\begin{equation}
c_{V,V}\left( e^i\otimes e^j\right) =R_{ab}^{ij}e^b\otimes e^a.  \label{x1}
\end{equation}

\noindent In the same way the category ${\cal M}^{H\left( R\right) }$ is
braided. If $W$ is the right fundamental corepresentation $W,$ with coaction 
$e_j\rightarrow e_a\otimes T_j^a,$ then the braiding is given by 
\begin{equation}
c_{W,W}\left( e_i\otimes e_j\right) =e_b\otimes e_aR_{ij}^{ab}.  \label{x2}
\end{equation}

\noindent If ${^{*}V}$ is the dual of the left fundamental corepresentation
with dual basis $\left\{ f_i\right\} $, then the coaction on ${^{*}V}$ is
given by $f_i\mapsto S\left( T_i^a\right) \otimes f_a.$ The braiding on
tensor products of $V$ and ${^{*}V}$ is now given by $\left( \ref{x1}\right) $
and the following formulas 
\begin{equation}
c_{{^{*}V},{^{*}V}}\left( f_i\otimes f_j\right) ={\em R}\left( ST_i^a\otimes
ST_j^b\right) f_b\otimes f_a=R_{ij}^{ab}f_b\otimes f_a,  \label{x3}
\end{equation}
\begin{equation}
c_{V,{^{*}V}}\left( e^i\otimes f_j\right) ={\em R}\left( T_a^i\otimes
ST_j^b\right) f_b\otimes e^a=\tilde R_{aj}^{ib}f_b\otimes e^a,  \label{x4}
\end{equation}
\begin{equation}
c_{{^{*}V},V}\left( f_i\otimes e^j\right) ={\em R}\left( ST_i^a\otimes
T_b^j\right) e^b\otimes f_a=\left( R^{-1}\right) _{ib}^{aj}e^b\otimes f_a.
\label{x5}
\end{equation}

\noindent Let ${^{*}W}$ be the dual of the right fundamental corepresentation,
with dual basis $\left\{ f^i\right\} $ and coaction $f^i\mapsto f^a\otimes
S\left( T_a^i\right) .$ Then the braiding is defined by $\left( \ref{x2}%
\right) $ and 
\[
c_{{^{*}W},{^{*}W}}\left( f^i\otimes f^j\right) ={\em R}\left( ST_a^i\otimes
ST_b^j\right) f^b\otimes f^a=f^b\otimes f^aR_{ab}^{ij}, 
\]
\[
c_{W,{^{*}W}}\left( e_i\otimes f^j\right) ={\em R}\left( T_i^a\otimes
ST_b^j\right) f^b\otimes e_a=f^b\otimes e_a\tilde R_{ib}^{aj}, 
\]
\[
c_{{^{*}W},W}\left( f^i\otimes e_j\right) ={\em R}\left( ST_b^i\otimes
T_j^a\right) e_a\otimes f^b=e_a\otimes f^b\left( R^{-1}\right) _{bj}^{ia}. 
\]
So $^{H(R)}{\cal M}$ and ${\cal M}^{H(R)}$ are braided categories with left
duality. The next definition dualizes the notion of ribbon algebra. 

\begin{definition}
Let $H$ be a dual quasi-triangular Hopf algebra with $2n^2$ generators $%
\left\{ x_1^1,\ldots ,x_n^n,y_1^1,\ldots ,y_n^n\right\} $ and 1, such that
it satisfies the hypotheses in corollary \ref{un}. Then we call $H$ a
coribbon algebra if $VU$ is invertible and has a central square root $\Theta
,$ which we call the coribbon element on $H,$ such that it defines an
element of $H^{*}$ by putting 
\[
\Theta \left( x_j^i\right) =\Theta _j^i,\,\,\Theta \left( y_j^i\right)
=\left( \Theta ^{-1}\right) _j^i,\,\,\Delta \left( \Theta \right) =\left( 
{\em R}_{21}{\em R}\right) ^{-1}\left( \Theta \otimes \Theta \right) , 
\]
\[
\,\varepsilon \left( \Theta \right) =1,\,\,S\left( \Theta \right) =\Theta . 
\]

\noindent Here we mean by the Hopf algebra structure the one in $H^{*}.$
\end{definition}

\noindent This gives us an analogue of theorem \ref{r}.

\begin{theorem}
Let $H$ satisfy the previous hypotheses. Then $H$ is a coribbon algebra {\it %
iff} $^H{\cal M}$ is a ribbon category ({\it iff} ${\cal M}^H$ is a ribbon
category).
\end{theorem}

\noindent {\bf Proof}. The proof is just the dual version of the proof of
theorem \ref{r}.\TeXButton{End Proof}{\hfill\endproof}

\medskip\ 

\noindent The next theorem dualizes theorem \ref{br}.

\begin{theorem}
\label{functor}let $R$ be a biinvertible R-matrix. If there is a coribbon
element $\Theta $ on $H(R)$, then $R$ and $\Theta $ define a unique functor $%
F_{RP}:{\cal FT}\rightarrow ^{H(R)}{\cal M}$ that preserves braiding,
duality and twist such that $F_{RP}\left( +\right) =V$ and $F_{RP}\left(
-\right) ={^{*}V}$ with $V$ the left fundamental corepresentation of $H(R).$
They also define a unique functor $F_{PR}:{\cal FT}\rightarrow {\cal M}%
^{H(R)}$ that preserves braiding, duality and twist such that $F_{PR}\left(
+\right) =W$ and $F_{PR}\left( -\right) =W^{*}$ with $W$ the right
fundamental corepresentation of $H(R).$
\end{theorem}

\noindent {\bf Proof}. Just dualize the proof of theorem \ref{br}.%
\TeXButton{End Proof}{\hfill\endproof}\ 

\medskip\ 

\noindent Note that we have to specify the coribbon element because there
might be more than one.

\begin{corollary}
\label{fun2} Let $R$ be a biinvertible R-matrix. If $\Theta =I$ defines a
coribbon element on $H(R),$ then there exists a unique functor $F_{RP}:{\cal %
T}\rightarrow ^{H(R)}{\cal M}$ that preserves braiding, duality and twist
such that $F_{RP}\left( +\right) =V$ and $F_{RP}\left( -\right) ={^{*}V}$ with 
$V$ the left fundamental corepresentation of $H(R).$ Furthermore there
exists a unique functor $F_{PR}:{\cal FT}\rightarrow {\cal M}^{H(R)}$ that
preserves braiding, duality and twist such that $F_{PR}\left( +\right) =W$
and $F_{PR}\left( -\right) ={^{*}W}$ with $W$ the right fundamental
corepresentation of $H(R).$
\end{corollary}

\noindent {\bf Proof.} Just dualize the proof of corollary \ref{univ2}.\hfill%
\TeXButton{End Proof}{\endproof}

\end{section} 

\begin{section}{Enhancement}

\noindent In this section we show exactly how enhancement fits into the
setup of the previous section. We first recall Turaev's theorem \cite{TUR3}. 

\begin{theorem}\label{Turaev}
Given an enhanced R-matrix in the sense of definition \ref{e2}%
, there exists a unique tensor functor $F\colon {\cal T}\rightarrow {\cal V_{%
{\Bbb C}}}$ such that $F\left( +\right) =V,F\left( -\right) =V^{*},$ and 
\[
\begin{tabular}{lll}
$F\left( X^{+}\right) =c,$ & $F\left( \cup \right) ={\rm coev}_V,$ & $%
F\left( \cup ^{-}\right) =\left( {\rm id}_{V^{*}}\otimes \mu ^{-1}\right) 
{\rm coev}_{V^{*}},$ \\ 
$F\left( X^{-}\right) =c^{-1},$ & $F\left( \cap \right) ={\rm ev}_V,$ & $%
F\left( \cap ^{-}\right) ={\rm ev}_{V^{*}}\left( \mu \otimes {\rm id}%
_{V^{*}}\right) .$%
\end{tabular}
\]
Conversely, let $F\colon {\cal T}\rightarrow {\cal V_{\Bbb C}}$ be a representation of ${\cal T}$ 
such that $F\left( +\right) =V, F\left( -\right)=V^{*}$ and 
\[
\begin{tabular}{lll} 
$F\left( X^{+}\right) =c,$ & $F\left( \cup \right) ={\rm coev}_V,$ & $%
F\left( \cup ^{-}\right) =b_{V}^{'},$ \\
$F\left( X^{-}\right) =c^{-1},$ & $F\left( \cap\right) ={\rm ev}_V,$ & $%
F\left( \cap ^{-}\right) =d_{V}^{'},$%
\end{tabular}
\]
where $c$ is an automorphism of $V\otimes V$. Then there exists a unique automorphism 
$\mu$ of $V$ such that $\left( c,\mu\right)$ is an enhanced R-matrix in the sense of def.\ref{e2} 
and 
$$b_{V}^{'}=\left( {\rm id}_{V^{*}}\otimes \mu^{-1}\right) {\rm coev}_{V^{*}},\ \ \ d_{V}^{'}={\rm ev}_{V^{*}}\left( \mu\otimes {\rm id}_{V^{*}}\right).$$
\end{theorem}
\noindent{\bf Sketch of a proof}. We only sketch the second part of 
the theorem, because we need it for the proof of our own results. For the rest of 
the proof we refer to \cite{TUR3} (see also \cite{KAS}). 
Let $\left\{ v_1,\ldots ,v_n\right\} $ be a basis of $%
V$ and $\left\{ w_1,\ldots ,w_n\right\} $ the dual basis of $V^{*}.$ We
define 
\[
b_V^{^{\prime }}\left( 1\right) =\sum_{i,j}B_{i,j}v_i\otimes w_j,%
\hspace{1.0in}d_V^{^{\prime }}\left( v_i\otimes w_j\right) =D_{i,j}. 
\]
\noindent Since $F$ is a representation of the tangle category we get 
\[
\left( d_V^{^{\prime }}\otimes {\rm id}_V\right) \left( {\rm id}_V\otimes
b_V^{^{\prime }}\right) ={\rm id}_V,\hspace{1.0in}\left( {\rm id}%
_{V^{*}}\otimes d_V^{^{\prime }}\right) \left( b_V^{^{\prime }}\otimes {\rm %
id}_{V^{*}}\right) ={\rm id}_{V^{*}}. 
\]
\noindent So $BD=DB=I.$ This allows us to define an automorphism $\beta
:V^{*}\rightarrow V^{*}$ by 
\[
\beta \left( w_j\right) =\sum_iB_{i,j}w_i. 
\]
\noindent Take $\mu =\left( \beta ^{-1}\right) ^t.$ Then reading Turaev's
proof shows that $\left( c,\mu \right) $ is the
enhanced R-matrix with the desired properties. The automorphism $\mu $ is
unique, because it is completely determined by the requirement $%
d_V^{^{\prime }}={\rm ev}_{V^{*}}\left( \mu \otimes {\rm id}_{V^{*}}\right)
. $\hfill\TeXButton{End Proof}{\endproof}

\noindent Next we give our results of this article.

\begin{theorem}
Let $R$ be a biinvertible R-matrix. If there exists a complex number $\alpha
\in {\bf {\Bbb C}}^{*}$ such that $VU=\alpha ^2I,$ then $\Theta =\alpha I$
is a coribbon element on $H(R).$ So $R$ and $\Theta $ define a unique
functor $F_{RP}:{\cal FT}\rightarrow ^{H(R)}{\cal M}$ that preserves
braiding, duality and twist such that $F_{RP}\left( +\right) =V$ and $%
F_{RP}\left( -\right) ={^{*}V}$ with $V$ the left fundamental corepresentation
of $H(R),$ and there is a unique functor $F_{PR}:{\cal FT}\rightarrow {\cal M%
}^{H(R)}$ that preserves braiding, duality and twist such that $F_{PR}\left(
+\right) =W$ and $F_{PR}\left( -\right) ={^{*}W}$ with $W$ the right
fundamental corepresentation of $H(R).$
\end{theorem}

\noindent {\bf Proof}. Obviously $\Theta =\alpha I$ defines a coribbon
element on $H(R)$ if we can prove that $\Theta $ is well defined in $H^{*}.$
We only prove $\Theta \left( TT^{-1}\right) =I,$ since the other identities
follow from the fact that $\Theta $ is central. 
\[
\Theta \left( T_\gamma ^i\left( T^{-1}\right) _l^\gamma \right) =\left(
\tilde R^{-1}\right) _{\alpha \beta }^{i\gamma }\left( \tilde R^{-1}\right)
_{ba}^{\beta \alpha }\Theta _\gamma ^a\Theta _l^b=\alpha ^2\left( \tilde
R^{-1}\right) _{\alpha \beta }^{i\gamma }\left( \tilde R\right) _{l\gamma
}^{\beta \alpha }. 
\]

\noindent Since we have 
\[
\left( \tilde R^{-1}\right) _{\alpha \beta }^{i\gamma }\left( \tilde
R^{-1}\right) _{l\gamma }^{\beta \alpha }\left( \tilde R\right) _{\gamma
\beta }^{l\gamma }\left( \tilde R\right) _{l\alpha }^{\alpha \beta }=\delta
_l^i 
\]

\noindent and 
\[
\left( \tilde R\right) _{\gamma \beta }^{l\gamma }\left( \tilde R\right)
_{l\alpha }^{\alpha \beta }=V_\beta ^lU_l^\beta =\left( VU\right)
_l^l=\alpha ^2, 
\]

\noindent we get 
\[
\Theta \left( T_\gamma ^i\left( T^{-1}\right) _l^\gamma \right) =\delta
_l^i. 
\]
So $\Theta $ defines a coribbon element on $H(R).$ The rest follows from
theorem \ref{functor}.\TeXButton{End Proof}{\hfill\endproof}

\medskip 

\begin{corollary}
\label{funct} Suppose $\alpha =1.$ Then there is a unique functor $F_{RP}:%
{\cal T}\rightarrow ^{H(R)}{\cal M}$ preserving braiding, duality and
twist such that $F_{RP}\left( +\right) =V$ and $F_{RP}\left( -\right) ={^{*}V}$
and a unique functor $F_{PR}:{\cal T}\rightarrow {\cal M}^{H(R)}$ that
preserves braiding, duality and twist such that $F_{PR}\left( +\right) =W$
and $F_{PR}\left( -\right) ={^{*}W}.$
\end{corollary}

\noindent {\bf Proof.} This follows from corollary \ref{fun2}.%
\TeXButton{End Proof}{\hfill\endproof}

\begin{corollary}
If $VU=\alpha ^2I$, then the matrix $R^{\prime }=\alpha R$ satisfies the
condition in the previous corollary.
\end{corollary}

\noindent {\bf Proof.} Trivial.

\medskip\ 

\noindent
\begin{theorem}
\label{mainres} If $R$ is a biinvertible R-matrix and $VU=I$, then $\left(
PR,U\right) $ and $\left( RP,V\right) $ are enhanced R-matrices in the sense
of definition \ref{e2}.
\end{theorem}

\noindent {\bf Proof.} The functors $F_{PR}$ and $F_{RP}$ in corollary \ref
{funct} satisfy the hypotheses in theorem \ref{Turaev}. In the first case we
get $\mu =\left( \beta ^{-1}\right) ^t=U$ and in the second case $\mu
=\left( \beta ^{-1}\right) ^t=V.$\hfill\TeXButton{End Proof}{\endproof}

\begin{corollary}
If $R$ is a biinvertible R-matrix and $VU=\alpha ^2I,$ then $\left( \alpha
PR,\alpha ^{-1}U\right) $ and $\left( \alpha RP,\alpha ^{-1}V\right) $ are
enhanced R-matrices in the sense of def.\ref{e2}. Of course $\left( PR,U,\alpha
^{-1},\alpha \right) $ and $\left( RP,V,\alpha ^{-1},\alpha \right) $ are
enhanced R-matrices in the sense of def.\ref{e1}.
\end{corollary}

\noindent {\bf Proof.} Obvious.\hfill\TeXButton{End Proof}{\endproof}

\medskip\ 

\noindent Our next lemma shows that the condition of $R$ being biinvertible
is no restriction at all, when one considers enhancement in the sense of def.%
\ref{e2}.

\begin{lemma}
\label{biinver} Let $R$ be an invertible R-matrix. If $R$ can be enhanced in
the sense of def.\ref{e2}, then $R$ is biinvertible.
\end{lemma}

\noindent {\bf Proof.} Suppose that $\left( PR,\mu \right) $ is enhanced in
the sense of def.\ref{e2}. One of the relations satisfied by the basic
elements of ${\cal T}$ (see \cite{KAS}) is the following 
\[
Y^{-}\circ T^{+}={\rm id}_{V^{*}}\otimes {\rm id}_V. 
\]

\noindent Note that this relation corresponds to an extended version of the
second Reidemeister move. When we apply the functor $F$ to the LHS of this
relation and use lemma \ref{rel} we find 
\[
F\left( Y^{-}\circ T^{+}\right) \left( f_a\otimes e^b\right) =R_{\beta
a}^{b\alpha }\left( \mu ^{-1}\right) _\delta ^\gamma \left( R^{-1}\right)
_{\varepsilon \theta }^{\delta \beta }\mu _\alpha ^\varepsilon f_\gamma
\otimes e^\theta . 
\]

\noindent If we take 
\[
A_{\theta \alpha }^{\beta \gamma }=\left( \mu ^{-1}\right) _\delta ^\gamma
\left( R^{-1}\right) _{\varepsilon \theta }^{\delta \beta }\mu _\alpha
^\varepsilon , 
\]

\noindent then we get 
\[
R_{\beta a}^{b\alpha }A_{\theta \alpha }^{\beta \gamma }=\delta _\gamma
^a\delta _\theta ^b. 
\]

\noindent Thus we see that $R$ is biinvertible with 
\[
A=\left( \left( R^{t_2}\right) ^{-1}\right) ^{t_2}=\tilde R. 
\]

\noindent There is a similar proof if $\left( RP,\nu \right) $ is enhanced.$%
\hfill\TeXButton{End Proof}{\endproof}$

\medskip\ 

\noindent Note that our proof also implies that the functor $F$ is actually a functor to the category of $H(R)$-comodules satisfying
\[
F\left( Y^{-}\right) =c_{V,V^{*}},\,\,F\left( T^{+}\right)
=c_{V^{*},V}^{-1}. 
\]

\noindent The same reasoning applied to the identity 
\[
Y^{+}\circ T^{-}={\rm id}_{V^{*}}\otimes {\rm id}_V 
\]

\noindent shows 
\[
F\left( Y^{+}\right) =c_{V,V^{*}}^{-1},\,\,F\left( T^{-}\right)
=c_{V^{*},V}. 
\]

\noindent Finally we can derive directly 
\[
F\left( Z^{\pm }\right) =c_{V^{*},V^{*}}^{\pm } 
\]
\noindent by using lemma \ref{rel}. Next we show that our condition is a
necessary condition for enhancement.

\begin{theorem}
\label{nec} If $\left( PR,\mu \right) $ is an enhanced R-matrix in the sense
of def.\ref{e2}, then $R$ is biinvertible and $VU=I.$ If $\left( RP,\nu
\right) $ is an enhanced R-matrix in the sense of def.\ref{e2}, then $R$ is
biinvertible and $VU=I.$
\end{theorem}

\noindent {\bf Proof}. Let us show the first case, the second being proven
in a similar way. If we take $V$ to be the fundamental left $H(R)$-comodule,
then theorem \ref{Turaev} shows that $\left( PR,\mu \right) $ defines a
functor $F$ from ${\cal T}$ to ${\cal V_{{\Bbb C}}}$ with $F\left( +\right)
=V$ and $F\left( -\right) =V^{*}.$ We define the subcategory $C(V)$ of $%
{\cal V_{{\Bbb C}}}$ generated by all tensor powers of $V$ and $V^{*}$. This
subcategory becomes a ribbon category if we define the braiding by the
formulas $\left( \ref{x1}\right) ,\left( \ref{x3}\right) ,\left( \ref{x4}%
\right) ,\left( \ref{x5}\right) $ and the left duality by the usual
evaluation and coevaluation maps. The twist is defined by $F\left( \varphi
\right) =$ $\theta _V:V\rightarrow V.$ We define $\theta _{^{*}V}$ and $%
\theta _{V\otimes V},\theta _{V\otimes V^{*}},\theta _{V^{*}\otimes
V},\theta _{V^{*}\otimes V^{*}}$ by the formulas \ref{t1} and \ref{t2}. Note 
that $C\left( V\right)$ is not a ribbon subcategory of ${\cal V}_{\Bbb C}$. For 
example $V^{**}\ne V$ as $H(R)$-comodules. The
previous lemma and the comments thereafter show that $R$ is biinvertible and
the functor $F$ maps ${\cal T}$ onto this subcategory, preserving braiding,
duality and twist. By lemma \ref{twist} we know how $\theta _V^{-2}$ acts.
Using the second expression in $\left( \ref{Theta}\right) $ we get 
\[
\begin{tabular}{ll}
$\theta _V^{-2}\left( e^i\right) $ & $=\left( {\rm ev}_Vc_{V,V^{*}}\otimes 
{\rm id}_V\right) \left( {\rm id}_V\otimes c_{V,V^{*}}{\rm coev}_V\right)
\left( e^i\right) $ \\ 
& $={\em R}\left( T_\alpha ^i\otimes ST_b^\beta \right) {\em R}\left(
T_a^j\otimes ST_j^b\right) \left\langle f_\beta ,e^\alpha \right\rangle
e^a=V_b^iU_a^be^a=\left( VU\right) _a^ie^a.$%
\end{tabular}
\]

\noindent So $\theta _V^{-2}=VU.$ Since Reidemeister move 1 is available in $%
{\cal T}$, we must have $\theta _V=I.$ Thus we get 
\[
VU=I. 
\]

\noindent By using the right fundamental $H(R)$-comodule and applying the
same arguments we get the result for the matrix $RP.$\hfill%
\TeXButton{End Proof}{\endproof}

\medskip\ 

\noindent Our last theorem proves that if a biinvertible R-matrix can be
enhanced, there is only one way to do it.

\begin{theorem}
\label{uniq} If $\left( PR,\mu \right) $ is an enhanced R-matrix in the
sense of def. \ref{e2}, then $R$ is biinvertible and $\mu =U.$ If $\left(
RP,\nu \right) $ is an enhanced R-matrix in the sense of def. \ref{e2}, then 
$R$ is biinvertible and $\nu =V.$
\end{theorem}

\noindent {\bf Proof.} Let $V$ be the left fundamental $H(R)$-comodule. If $%
\left( PR,\mu \right) $ is an enhanced R-matrix, then, as we showed in the
previous theorem, $R$ is biinvertible and there is a unique functor $F:{\cal %
T}\rightarrow C(V)$ that preserves braiding, left duality and twist. Of
course $F$ also preserves right duality, so 
\begin{equation}
F\left( \cup ^{-}\right) =\left( {\rm id}_{V^{*}}\otimes \mu ^{-1}\right) 
{\rm coev}_{V^{*}},\,\,\,F\left( \cap ^{-}\right) ={\rm ev}_{V^{*}}\left(
\mu \otimes {\rm id}_{V^{*}}\right)
\end{equation}

\noindent define a right duality in $C(V).$ On the other hand we know that 
\[
Y^{-}\circ \cup =\cup ^{-},\,\,\cap \circ Y^{-}=\cap ^{-} 
\]
\noindent are satisfied in ${\cal T}$. Hence 
\[
F\left( \cup ^{-}\right) =F\left( Y^{-}\right) F\left( \cup \right)
=c_{V,V^{*}}{\rm coev}_V,\,\,F\left( \cap ^{-}\right) =F\left( \cap \right)
F\left( Y^{-}\right) ={\rm ev}_Vc_{V,V^{*}}. 
\]
\noindent So both right dualities are the same; 
\[
\left( {\rm id}_{V^{*}}\otimes \mu ^{-1}\right) {\rm coev}%
_{V^{*}}=c_{V,V^{*}}{\rm coev}_V, 
\]
\[
{\rm ev}_{V^{*}}\left( \mu \otimes {\rm id}_{V^{*}}\right) ={\rm ev}%
_Vc_{V,V\ ^{*}}. 
\]

\noindent As a result we get 
\[
\mu _j^i={\rm ev}_{V^{*}}\left( \mu \otimes {\rm id}_{V^{*}}\right) \left(
e^i\otimes f_j\right) ={\rm ev}_Vc_{V,V^{*}}\left( e^i\otimes f_j\right)
=\tilde R_{aj}^{ib}\left\langle f_{b,}e^a\right\rangle =U_j^i. 
\]

\noindent In the same way one can prove the statement for the matrix $RP.$%
\hfill\TeXButton{End Proof}{\endproof}

\medskip\ 

\noindent Note that theorem \ref{nec} and \ref{uniq} are only true for
enhancement in the sense of def.\ref{e2}. If $R$ is biinvertible and $\left(
PR,U,1,1\right) $ is enhanced in the sense of def.\ref{e1}, then $UV={\rm Tr}%
_2\left( \left( PR\right) ^{-1}U_2\right) =I$ (see Appendix). If $R$ is
biinvertible and $\left( RP,V,1,1\right) $ is enhanced in the sense of def.%
\ref{e1}, then $VU={\rm Tr}_2\left( \left( RP\right) ^{-1}V_2\right) =I$
(see Appendix). But we don't know if there are no other $\mu $ and $\nu $ such
that $\left( PR,\mu ,1,1\right) $ and $\left( RP,\nu ,1,1\right) $ are
enhanced.

\bigskip\ 

\noindent 5. {\bf Examples}

\medskip\ 

\noindent As an example we enhance all biinvertible R-matrices of dimension
2.

\[
R=\left( 
\begin{array}{cccc}
1 & 0 & 0 & 0 \\ 
0 & p & 0 & 0 \\ 
0 & 0 & s & 0 \\ 
0 & 0 & 0 & q
\end{array}
\right) ,\,\tilde R=\left( 
\begin{array}{cccc}
1 & 0 & 0 & 0 \\ 
0 & p^{-1} & 0 & 0 \\ 
0 & 0 & s^{-1} & 0 \\ 
0 & 0 & 0 & q^{-1}
\end{array}
\right) , 
\]

\[
U=\left( 
\begin{array}{cc}
1 & 0 \\ 
0 & q^{-1}
\end{array}
\right) ,\,V=\left( 
\begin{array}{cc}
1 & 0 \\ 
0 & q^{-1}
\end{array}
\right) . 
\]

\noindent For $q=1$ the matrices $\left( PR,U\right) $ and $\left(
RP,V\right) $ are enhanced. 
\[
R=\left( 
\begin{array}{cccc}
0 & 0 & 0 & q \\ 
0 & 0 & 1 & 0 \\ 
0 & 1 & 0 & 0 \\ 
q & 0 & 0 & 0
\end{array}
\right) ,\,\tilde R=\left( 
\begin{array}{cccc}
0 & 0 & 0 & q^{-1} \\ 
0 & 0 & 1 & 0 \\ 
0 & 1 & 0 & 0 \\ 
q^{-1} & 0 & 0 & 0
\end{array}
\right) 
\]
$,$%
\[
U=V=\left( 
\begin{array}{cc}
1 & 0 \\ 
0 & 1
\end{array}
\right) . 
\]

\noindent The matrices $\left( PR,U\right) $ and $\left( RP,V\right) $ are
enhanced. 
\[
R=\left( 
\begin{array}{cccc}
1 & 1 & p & q \\ 
0 & 1 & 0 & p \\ 
0 & 0 & 1 & 1 \\ 
0 & 0 & 0 & 1
\end{array}
\right) ,\,\tilde R=\left( 
\begin{array}{cccc}
1 & -1 & -p & 2p-q \\ 
0 & 1 & 0 & -p \\ 
0 & 0 & 1 & -1 \\ 
0 & 0 & 0 & 1
\end{array}
\right) , 
\]
\[
U=V=\left( 
\begin{array}{cc}
1 & -p-1 \\ 
0 & 1
\end{array}
\right) . 
\]

\noindent For $p=-1$ the matrices $\left( PR,U\right) $ and $\left(
RP,V\right) $ are enhanced. 
\[
R=\left( 
\begin{array}{cccc}
1 & 1 & -1 & q \\ 
0 & 1 & 0 & q \\ 
0 & 0 & 1 & -q \\ 
0 & 0 & 0 & 1
\end{array}
\right) ,\,\tilde R=\left( 
\begin{array}{cccc}
1 & -1 & 1 & -q^2-q-1 \\ 
0 & 1 & 0 & -q \\ 
0 & 0 & 1 & q \\ 
0 & 0 & 0 & 1
\end{array}
\right) , 
\]
\[
U=\left( 
\begin{array}{cc}
1 & 1+q \\ 
0 & 1
\end{array}
\right) ,\,V=\left( 
\begin{array}{cc}
1 & -1-q \\ 
0 & 1
\end{array}
\right) . 
\]

\noindent $\left( PR,U\right) $ and $\left( RP,V\right) $ are enhanced. 
\[
R=\left( 
\begin{array}{cccc}
1 & 0 & 0 & 1 \\ 
0 & 1 & 1 & 0 \\ 
0 & 1 & -1 & 0 \\ 
-1 & 0 & 0 & 1
\end{array}
\right) 
\]

\noindent $R$ is not biinvertible. 
\[
R=\left( 
\begin{array}{cccc}
1 & 0 & 0 & 1 \\ 
0 & 1 & 0 & 0 \\ 
0 & 0 & 1 & 0 \\ 
0 & 0 & 0 & 1
\end{array}
\right) ,\,\tilde R=\left( 
\begin{array}{cccc}
1 & 0 & 0 & -1 \\ 
0 & 1 & 0 & 0 \\ 
0 & 0 & 1 & 0 \\ 
0 & 0 & 0 & 1
\end{array}
\right) , 
\]
\[
U=V=\left( 
\begin{array}{cc}
1 & 0 \\ 
0 & 1
\end{array}
\right) . 
\]

\noindent $\left( PR,U\right) $ and $\left( RP,V\right) $ are enhanced. 
\[
R=\left( 
\begin{array}{cccc}
1 & 0 & 0 & 1 \\ 
0 & -1 & 0 & 0 \\ 
0 & 0 & -1 & 0 \\ 
0 & 0 & 0 & 1
\end{array}
\right) ,\,\tilde R=\left( 
\begin{array}{cccc}
1 & 0 & 0 & -1 \\ 
0 & -1 & 0 & 0 \\ 
0 & 0 & -1 & 0 \\ 
0 & 0 & 0 & 1
\end{array}
\right) , 
\]
\[
U=V=\left( 
\begin{array}{cc}
1 & 0 \\ 
0 & 1
\end{array}
\right) . 
\]

\noindent $\left( PR,U\right) $ and $\left( RP,V\right) $ are enhanced. 
\[
R=\left( 
\begin{array}{cccc}
q & 0 & 0 & 0 \\ 
0 & p & q-q^{-1} & 0 \\ 
0 & 0 & p^{-1} & 0 \\ 
0 & 0 & 0 & q
\end{array}
\right) ,\,\tilde R=\left( 
\begin{array}{cccc}
q^{-1} & 0 & 0 & 0 \\ 
0 & p^{-1} & \frac{q^{-1}-q}{q^2} & 0 \\ 
0 & 0 & p & 0 \\ 
0 & 0 & 0 & q^{-1}
\end{array}
\right) , 
\]
\[
U=\left( 
\begin{array}{cc}
q^{-1} & 0 \\ 
0 & q^{-3}
\end{array}
\right) ,\,V=\left( 
\begin{array}{cc}
q^{-3} & 0 \\ 
0 & q^{-1}
\end{array}
\right) . 
\]

\noindent $\left( q^{-2}PR,q^2U\right) $ and $\left( q^{-2}RP,q^2V\right) $
are enhanced. 

\[
R=\left( 
\begin{array}{cccc}
q & 0 & 0 & q \\ 
0 & p & q-q^{-1} & 0 \\ 
0 & 0 & p^{-1} & 0 \\ 
0 & 0 & 0 & -q^{-1}
\end{array}
\right) ,\,\tilde R=\left( 
\begin{array}{cccc}
q^{-1} & 0 & 0 & -q \\ 
0 & p^{-1} & q-q^{-1} & 0 \\ 
0 & 0 & p & 0 \\ 
0 & 0 & 0 & -q
\end{array}
\right) 
\]
\[
U=\left( 
\begin{array}{cc}
q^{-1} & 0 \\ 
0 & -q^{-1}
\end{array}
\right) ,\,V=\left( 
\begin{array}{cc}
q & 0 \\ 
0 & -q
\end{array}
\right) . 
\]

\noindent $\left( PR,U\right) $ and $\left( RP,V\right) $ are enhanced. 
\[
R=\left( 
\begin{array}{cccc}
q-q^{-1}+2 & 0 & 0 & q-q^{-1} \\ 
0 & q+q^{-1} & q-q^{-1} & 0 \\ 
0 & q-q^{-1} & q+q^{-1} & 0 \\ 
q-q^{-1} & 0 & 0 & q-q^{-1}-2
\end{array}
\right) , 
\]
\[
\tilde R=\frac 14\left( 
\begin{array}{cccc}
q-q^{-1}+2 & 0 & 0 & q^{-1}-q \\ 
0 & q+q^{-1} & q-q^{-1} & 0 \\ 
0 & q-q^{-1} & q+q^{-1} & 0 \\ 
q^{-1}-q & 0 & 0 & q^{-1}-q-2
\end{array}
\right) , 
\]
\[
U=V=\frac 12\left( 
\begin{array}{cc}
1 & 0 \\ 
0 & -1
\end{array}
\right) . 
\]

\noindent $\left( \frac 12PR,2U\right) $ and $\left( \frac 12RP,2V\right) $
are enhanced.

\noindent The permutation matrix 
\[
P=\left( 
\begin{array}{cccc}
1 & 0 & 0 & 0 \\ 
0 & 0 & 1 & 0 \\ 
0 & 1 & 0 & 0 \\ 
0 & 0 & 0 & 1
\end{array}
\right) 
\]

\noindent is not biinvertible, as one can check easily. But, as we wrote in
the introduction already, $\left( P^2,I,\frac 12,2\right) $ is an enhanced
matrix.

\bigskip\ 

\noindent {\bf Appendix}

\medskip\ 

\noindent In this appendix we give an elementary proof of corollary \ref
{mainres}.

\begin{theorem}
Let $R$ be a biinvertible R-matrix. If $VU=UV=I$, then $\left( RP,V\right) $
and $\left( PR,U\right) $ are enhanced R-matrices in the sense of def.\ref
{e2}.
\end{theorem}

\noindent {\bf Proof}. In order to prove that $\left( RP,V\right) $ is an
enhanced R-matrix we use the identities $V_1=\tilde RV_1R$ and $%
V_2=RV_2\tilde R$ $\left( \text{see \cite{MAJ}}\right) $. We first prove
that $V\otimes V=V_1V_2=V_2V_1$ commutes with $R$: 
\[
RPV_1V_2=RV_2PV_2=RV_2V_1P=V_2\tilde R_1^{-1}VP=V_2V_1RP. 
\]

\noindent Next we prove the conditions on the partial traces: 
\[
\begin{array}{l}
\left( \text{Tr}_2\left( RPV_2\right) \right) _c^a=\left( RPV_2\right)
_{cd}^{ad}=R_{ji}^{ad}\left( {\rm id}\otimes V\right)
_{cd}^{ij}=R_{ji}^{ad}\delta _c^iV_d^j=R_{jc}^{ad}\tilde R_{kd}^{jk}=\delta
_k^a\delta _c^k=\delta _c^a, \\ 
\\ 
\left( \text{Tr}_2\left( \left( RP\right) ^{-1}V_2\right) \right)
_c^a=\left( \left( RP\right) ^{-1}V_2\right) _{cd}^{ad}=\left(
PR^{-1}V_2\right) _{cd}^{ad}=\left( PR^{-1}RV_2\tilde R\right) _{cd}^{ad} \\ 
=\left( PV_2\tilde R\right) _{cd}^{ad}=\left( {\rm id}\otimes V\right)
_{ij}^{da}\tilde R_{cd}^{ij}=\delta _i^dV_j^a\tilde R_{cd}^{ij}=\tilde
R_{kj}^{ak}\tilde R_{cd}^{dj}=V_j^aU_c^j=\delta _c^a.
\end{array}
\]

\noindent Finally we prove condition$\left( \ref{ENH4}\right) $. It's not difficult (see \cite{TUR3} or \cite{KAS}) to show that this 
condition is equivalent to the following condition 
\begin{equation}  \label{ENH5}
\left( {\rm Id}_V\otimes \left( \mu ^{*}\right) ^{-1}\right) \left( S^{\pm
1}P\right) ^{t_2}\left( {\rm id}_V\otimes \mu ^{*}\right) \left( PS^{\pm
1}\right) ^{t_2}={\rm id}_{V\otimes {^{*}V}}.
\end{equation}

\noindent We prove that $\left( RP,V\right) $ satisfies $\left( \ref{ENH5}%
\right) $. 
\begin{eqnarray*}
\left( V_2^{t_2}\right) ^{-1}R^{t_2}V_2^{t_2}\left( R^{-1}\right) ^{t_2}
&=&\left( V_2^{t_2}\right) ^{-1}R^{t_2}\left( R^{-1}V_2\right) ^{t_2}= \\
\left( V_2^{t_2}\right) ^{-1}R^{t_2}\left( V_2\tilde R\right) ^{t_2}
&=&\left( V_2^{t_2}\right) ^{-1}R^{t_2}\left( R^{t_2}\right) ^{-1}V_2^{t_2}=I
\end{eqnarray*}

\noindent This proves the first identity in $\left( \ref{ENH5}\right) $. The
second follows from 
\[
\begin{array}{c}
\left( PR^{-1}P\right) ^{t_2}V_2^{t_2}=\left( V_2PR^{-1}P\right)
^{t_2}=\left( PV_1R^{-1}P\right) ^{t_2}=\left( P\tilde RV_1P\right) ^{t_2}
\\ 
=\left( P\tilde RPV_2\right) ^{t_2}=V_2^{t_2}\left( P\tilde RP\right)
^{t_2}=V_2^{t_2}\left( \left( PRP\right) ^{t_2}\right) ^{-1}.
\end{array}
\]

\noindent The other pair is proven to be an enhanced R-matrix in the
same way, using the identities $U_1=RU_1\tilde R$ and $U_2=\tilde RU_2R.$ 
\[
PRU_1U_2=PU_1\tilde R_2^{-1}U_2=PU_1U_2R=U_2PU_2R=U_2U_1PR 
\]
\[
\begin{array}{l}
\left( \left( \text{Tr}_2\left( PRU_2\right) \right) ^t\right) _c^a=\text{Tr}%
_2\left( \left( PRU_2\right) ^t\right) _c^a=\text{Tr}_2\left( U_2^t\left(
PR\right) ^t\right) _c^a=\left( {\rm id}\otimes U\right) _{ad}^{ij}\left(
PR\right) _{ij}^{cd} \\ 
=\delta _a^i\tilde R_{dk}^{kj}R_{ij}^{dc}=\tilde
R_{dk}^{kj}R_{aj}^{dc}=\delta _a^k\delta _k^c=\delta _a^c \\ 
\\ 
\left( \text{Tr}_2\left( \left( PR\right) _2^{-1}U_2\right) \right)
_c^a=\left( \left( PR\right) _2^{-1}U_2\right) _{cd}^{ad}=\left( \left(
R^{-1}PU_2\right) \right) _{cd}^{ad}=\left( R^{-1}U_1P\right) _{cd}^{ad} \\ 
=\left( U_1\tilde RP\right) _{cd}^{ad}=(U\otimes {\rm id})_{ij}^{ad}\tilde
R_{dc}^{ij}=U_i^a\delta _j^d\tilde R_{dc}^{ij}=U_i^aV_c^i=\delta _c^a.
\end{array}
\]

\noindent The proof that $\left( PR,U\right) $ satisfies $\left( \ref{ENH5}%
\right) $ follows in an analogous way from the identities $U_1=RU_1\tilde R$
and $U_2=\tilde RU_2R.$\TeXButton{End Proof}{\hfill\endproof}

\medskip\ 

\begin{corollary}
\label{main}Let $R$ be a biinvertible R-matrix. If there exists a scalar $%
\alpha \in {\bf {\Bbb C}}^{*}$ such that 
\[
\alpha ^2UV=I, 
\]

\noindent then $\left( \alpha PR,\alpha ^{-1}U\right) $ and $\left( \alpha
RP,\alpha ^{-1}V\right) $ are enhanced R-matrices in the sense of def.\ref
{e2} and $\left( PR,U,\alpha ^{-1},\alpha \right) $ and $\left( RP,V,\alpha
^{-1},\alpha \right) $ are enhanced in the sense of def.\ref{e1}.
\end{corollary}

\noindent {\bf Proof}. This clearly follows from the theorem above.%
\TeXButton{End Proof}{\hfill\endproof}
\end{section}

\end{document}